\documentclass[11pt]{amsart}
\usepackage{amssymb,amsmath}
\usepackage{graphicx}
\usepackage[all]{xy}

\title[Invariance of the underlying homotopy type]{T-homotopy and refinement of observation (IV) : Invariance of the underlying homotopy type}
\author[P. Gaucher]{Philippe Gaucher}
\address{Preuves Programmes et Syst{\`e}mes\\ Universit{\'e} Paris 7--Denis Diderot\\
Case 7014\\2 Place Jussieu\\ 75251 PARIS Cedex 05\\ France}
\email{gaucher@pps.jussieu.fr}
\urladdr{http://www.pps.jussieu.fr/{\~{}}gaucher/}

\newcommand{\C}{\mathcal{C}}
\newcommand{\N}{\mathbb{N}}
\newcommand{\de}{\partial}
\newcommand{\p}\times
\renewcommand{\vec}{\overrightarrow}
\renewcommand{\P}{\mathbb{P}}

\newcommand{\be}{\begin{equation}}
\newcommand{\ee}{\end{equation}}
\newcommand{\beas}{\begin{eqnarray*}}
\newcommand{\eeas}{\end{eqnarray*}}

\newtheorem*{thmN}{Theorem}
\newtheorem{thm}{Theorem}[section]
\newtheorem{prop}[thm]{Proposition}
\newtheorem{lem}[thm]{Lemma}
\newtheorem{cor}[thm]{Corollary}
\newtheorem{defn}[thm]{Definition}
\newtheorem{nota}[thm]{Notation}
\newcommand{\bd}{\begin{defn}}
\newcommand{\ed}{\end{defn}}
\newcommand{\bp}{\begin{prop}}
\newcommand{\ep}{\end{prop}}
\newcommand{\bth}{\begin{thm}}
\renewcommand{\eth}{\end{thm}}
\newcommand{\bpf}{\begin{proof}}
\newcommand{\epf}{\end{proof}}

\newcommand{\fl}[1]{\ar@{->}[l]_{#1}}
\newcommand{\fr}[1]{\ar@{->}[r]^-{#1}}
\newcommand{\fd}[1]{\ar@{->}[d]_{#1}}
\newcommand{\fu}[1]{\ar@{->}[u]^{#1}}

\newcommand{\ho}{{\mathbf{Ho}}}

\newcommand{\iso}{\cong}
\newcommand{\lp}{\left(}
\newcommand{\rp}{\right)}

\newcommand{\vI}{\vec{I}}

\renewcommand{\leq}{\leqslant}
\renewcommand{\geq}{\geqslant}

\def\cartesien{%
  \ar@{-}[]+R+<6pt,-2pt>;[]+RD+<6pt,-6pt>%
  \ar@{-}[]+D+<2pt,-6pt>;[]+RD+<6pt,-6pt>%
}
\def\cocartesien{%
  \ar@{-}[]+L+<-6pt,+2pt>;[]+LU+<-6pt,+6pt>%
  \ar@{-}[]+U+<-2pt,+6pt>;[]+LU+<-6pt,+6pt>%
}

\newcommand{\brm}[1]{\rm{\mathbf{#1}}}

\renewcommand{\top}{{\brm{Top}}}
\newcommand{\gltop}{{\brm{glTop}}}
\newcommand{\dtop}{{\brm{Flow}}}

\newcommand{\set}{{\brm{Set}}}
\newcommand{\glob}{{\rm{Glob}}}
\newcommand{\liminj}{\varinjlim}

\makeatletter

\def\varholim@#1#2{%
  \vtop{\m@th\ialign{##\cr
    \hfil$#1\operator@font holim$\hfil\cr
    \noalign{\nointerlineskip\kern1.5\ex@}#2\cr
    \noalign{\nointerlineskip\kern-\ex@}\cr}}%
}
\def\holimproj{%
  \mathop{\mathpalette\varholim@{\leftarrowfill@\textstyle}}\nmlimits@
}
\def\holiminj{%
  \mathop{\mathpalette\varholim@{\rightarrowfill@\textstyle}}\nmlimits@
}

\makeatother

\DeclareMathOperator{\id}{Id}
\DeclareMathOperator{\card}{card}

\DeclareMathOperator{\diag}{{\rm{Diag}}}

\DeclareMathOperator{\cell}{{\brm{cell}}}
\DeclareMathOperator{\cof}{{\brm{cof}}}
\DeclareMathOperator{\inj}{{\brm{inj}}}
\newcommand{\hda}{{\cell(\dtop)}}

\begin{document}

\begin{abstract} 
  This series explores a new notion of T-homotopy equivalence of
  flows.  The new definition involves embeddings of finite bounded
  posets preserving the bottom and the top elements and the associated
  cofibrations of flows. In this fourth part, it is proved that the
  generalized T-homotopy equivalences preserve the underlying homotopy
  type of a flow. The proof is based on Reedy model category
  techniques.
\end{abstract}

\maketitle

\tableofcontents

\section{Outline of the paper}

The main feature of the two algebraic topological models of
\textit{higher dimensional automata} (or HDA) introduced in
\cite{diCW} and in \cite{model3} is to provide a framework for
modelling continuous deformations of HDA corresponding to subdivision
or refinement of observation. \textit{Globular complexes} and
\textit{flows} are introduced in \cite{diCW} and \cite{model3}
respectively for modelling a notion of dihomotopy equivalence between
higher dimensional automata \cite{Pratt} \cite{rvg}.  This equivalence
relation preserves geometric properties like the \textit{initial} or
\textit{final states}, and therefore computer-scientific properties
like the presence or not of \textit{deadlocks} or of
\textit{unreachable states} \cite{survol}.  More generally, dihomotopy
is designed to preserving all computer-scientific properties invariant
by refinement of observation (see Figure~\ref{ex1}). The two settings
are compared in \cite{model2} and are proved to be equivalent.

In the framework of flows, there are two kinds of dihomotopy
equivalences \cite{ConcuToAlgTopo}: the \textit{weak S-homotopy
  equivalences} (the spatial deformations of \cite{ConcuToAlgTopo})
which can be interpreted as the weak equivalences of a model structure
\cite{model3} and the \textit{T-homotopy equivalences} (the temporal
deformations of \cite{ConcuToAlgTopo}). The latter are considerably
more difficult to model and to understand. The geometric explanations
underlying the intuition of S-homotopy and T-homotopy are given in the
first part of this series \cite{1eme}, but the reference \cite{diCW}
must be preferred.

The purpose of this paper is to prove that the notion of T-homotopy
equivalence studied in this series preserves the \textit{underlying
  homotopy type} of a flow. The underlying homotopy type of a flow is
the topological space which is obtained after removing the temporal
ordering. This underlying topological space is unique only up to weak
homotopy equivalence. For example, the underlying homotopy type of the
two flows of Figure~\ref{ex1} is the point. The main theorem of this
paper is:
\begin{thmN} 
Let $f:X\longrightarrow Y$ be a generalized T-homotopy equivalence.
Then the morphism of $\ho(\top)$ $|f|:|X|\longrightarrow |Y|$, where
$|-|$ is the underlying homotopy type functor, is an isomorphism.
\end{thmN}

Section~\ref{remindergg} recalls the notions of full directed ball and
of generalized T-homotopy equivalence. Section~\ref{remindergl}
recalls the notion of globular complex. It is necessary for the
definition of the underlying homotopy type of a flow.
Section~\ref{defunder} gives the rigorous definition of the underlying
homotopy type of a flow. Section~\ref{reedyuse} constructs a useful
Reedy structure which will be crucial in the main proofs of the paper.
Section~\ref{reedyuse} also establishes related lemmas.
Section~\ref{calhoty} proves that the underlying homotopy type of the
full directed ball is contractible (Theorem~\ref{fin}).  The latter
result is important since a T-homotopy equivalence consists in
replacing in a flow a full directed ball by a more refined full
directed ball (see Figure~\ref{ex2}), and in iterating this process
transfinitely. Then Section~\ref{prehoty} proves the theorem above.

\subsection*{Warning.} This paper is the fourth part of a
series of papers devoted to the study of T-homotopy. Several other
papers explain the geometrical content of T-homotopy. The best
reference is probably \cite{diCW} (it does not belong to the series).
However, the knowledge of the other parts is not required. In
particular, this means that there are repetitions between the papers
of this series. They are all of them collected in the appendices
\ref{limelm}, \ref{calpush} and \ref{mixcomp} which are already in the
third part of this series. The proofs of these appendices are
independent from the technical core of this part.  The left properness
of the weak S-homotopy model structure of $\dtop$ is not duplicated in
this paper. It is available in \cite{2eme}. This fact is used twice in
the proof of Theorem~\ref{thth}.

\section{Prerequisites and notations}

The initial object (resp. the terminal object) of a category $\C$, if
it exists, is denoted by $\varnothing$ (resp. $\mathbf{1}$).

Let $\C$ be a cocomplete category.  If $K$ is a set of morphisms of
$\C$, then the class of morphisms of $\C$ that satisfy the RLP
(\textit{right lifting property}) with respect to any morphism of $K$
is denoted by $\inj(K)$ and the class of morphisms of $\C$ that are
transfinite compositions of pushouts of elements of $K$ is denoted by
$\cell(K)$. Denote by $\cof(K)$ the class of morphisms of $\C$ that
satisfy the LLP (\textit{left lifting property}) with respect to the
morphisms of $\inj(K)$.  This is a purely categorical fact that
$\cell(K)\subset \cof(K)$. Moreover, every morphism of $\cof(K)$ is a
retract of a morphism of $\cell(K)$ as soon as the domains of $K$ are
small relative to $\cell(K)$ (\cite{MR99h:55031} Corollary~2.1.15). An
element of $\cell(K)$ is called a \textit{relative $K$-cell complex}.
If $X$ is an object of $\C$, and if the canonical morphism
$\varnothing\longrightarrow X$ is a relative $K$-cell complex, then
the object $X$ is called a \textit{$K$-cell complex}.

Let $\C$ be a cocomplete category with a distinguished set of
morphisms $I$. Then let $\cell(\C,I)$ be the full subcategory of $\C$
consisting of the objects $X$ of $\C$ such that the canonical morphism
$\varnothing\longrightarrow X$ is an object of $\cell(I)$. In other
terms, $\cell(\C,I)=(\varnothing\!\downarrow \! \C) \cap \cell(I)$.

It is obviously impossible to read this paper without a strong
familiarity with \textit{model categories}. Possible references for
model categories are \cite{MR99h:55031}, \cite{ref_model2} and
\cite{MR1361887}.  The original reference is \cite{MR36:6480} but
Quillen's axiomatization is not used in this paper. The axiomatization
from Hovey's book is preferred.  If $\mathcal{M}$ is a
\textit{cofibrantly generated} model category with set of generating
cofibrations $I$, let $\cell(\mathcal{M}) := \cell(\mathcal{M},I)$ :
this is the full subcategory of \textit{cell complexes} of the model
category $\mathcal{M}$. A cofibrantly generated model structure
$\mathcal{M}$ comes with a \textit{cofibrant replacement functor}
$Q:\mathcal{M} \longrightarrow \cell(\mathcal{M})$. For any morphism
$f$ of $\mathcal{M}$, the morphism $Q(f)$ is a cofibration, and even
an inclusion of subcomplexes (\cite{ref_model2} Definition~10.6.7)
because the cofibrant replacement functor $Q$ is obtained by the small
object argument.

A \textit{partially ordered set} $(P,\leq)$ (or \textit{poset}) is a
set equipped with a reflexive antisymmetric and transitive binary
relation $\leq$. A poset is \textit{locally finite} if for any
$(x,y)\in P\p P$, the set $[x,y]=\{z\in P,x\leq z\leq y\}$ is finite.
A poset $(P,\leq)$ is \textit{bounded} if there exist $\widehat{0}\in
P$ and $\widehat{1}\in P$ such that $P = [\widehat{0},\widehat{1}]$
and such that $\widehat{0} \neq \widehat{1}$. Let $\widehat{0}=\min P$
(the bottom element) and $\widehat{1}=\max P$ (the top element). In a
poset $P$, the interval $]\alpha,-]$ (the sub-poset of elements of $P$
strictly bigger than $\alpha$) can also be denoted by $P_{>\alpha}$.

A poset $P$, and in particular an ordinal, can be viewed as a small
category denoted in the same way: the objects are the elements of $P$
and there exists a morphism from $x$ to $y$ if and only if $x\leq
y$. If $\lambda$ is an ordinal, a \textit{$\lambda$-sequence} in a
cocomplete category $\C$ is a colimit-preserving functor $X$ from
$\lambda$ to $\C$. We denote by $X_\lambda$ the colimit $\liminj X$
and the morphism $X_0\longrightarrow X_\lambda$ is called the
\textit{transfinite composition} of the $X_\mu\longrightarrow
X_{\mu+1}$.

Let $\C$ be a category. Let $\alpha$ be an object of $\C$. The
\textit{latching category} $\de(\C\!\downarrow\! \alpha)$ at $\alpha$ is the
full subcategory of $\C\!\downarrow\! \alpha$ containing all the
objects except the identity map of $\alpha$. The
\textit{matching category} $\de(\alpha\!\downarrow\!\C)$ at $\alpha$ is the
full subcategory of $\alpha\!\downarrow\!\C$ containing all the
objects except the identity map of $\alpha$.

Let $\mathcal{B}$ be a small category. A \textit{Reedy structure} on
$\mathcal{B}$ consists of two subcategories $\mathcal{B}_-$ and
$\mathcal{B}_+$, a functor $d:\mathcal{B}\longrightarrow \lambda$
called the {\rm degree function} for some ordinal $\lambda$, such that
every non identity map in $\mathcal{B}_+$ raises the degree, every
non identity map in $\mathcal{B}_-$ lowers the degree, and every map
$f\in \mathcal{B}$ can be factored uniquely as $f=g \circ h$ with $h
\in \mathcal{B}_-$ and $g \in \mathcal{B}_+$. A small category together
with a Reedy structure is called a \textit{Reedy category}.

Let $\C$ be a complete and cocomplete category. Let $\mathcal{B}$ be a
Reedy category. Let $i$ be an object of $\mathcal{B}$. The \textit{latching
space functor} is the composite $L_i:\C^\mathcal{B}\longrightarrow
\C^{\de(\mathcal{B}_+\!\downarrow\! i)}\longrightarrow \C$ where the latter
functor is the colimit functor.  The \textit{matching space functor}
is the composite $M_i:\C^\mathcal{B}\longrightarrow
\C^{\de(i\!\downarrow\!\mathcal{B}_-)}\longrightarrow \C$ where the latter
functor is the limit functor.

If $\C$ is a small category and of $\mathcal{M}$ is a model category,
the notation $\mathcal{M}^\C$ is the category of functors from $\C$ to
$\mathcal{M}$, i.e. the category of diagrams of objects of
$\mathcal{M}$ over the small category $\C$.

The category $\top$ of \textit{compactly generated topological spaces}
(i.e. of weak Hausdorff $k$-spaces) is complete, cocomplete and
cartesian closed (more details for this kind of topological spaces in
\cite{MR90k:54001,MR2000h:55002}, the appendix of \cite{Ref_wH} and
also the preliminaries of \cite{model3}). For the sequel, all
topological spaces will be supposed to be compactly generated. A
\textit{compact space} is always Hausdorff.

A model category is \textit{left proper} if the pushout of a weak
equivalence along a cofibration is a weak equivalence. The model
categories $\top$ and $\dtop$ (see below) are both left proper.

In this paper, the notation $\xymatrix@1{\ar@{^{(}->}[r]&}$ means
\textit{cofibration}, the notation $\xymatrix@1{\ar@{->>}[r]&}$ means
\textit{fibration}, the notation $\simeq$ means \textit{weak
  equivalence}, and the notation $\iso$ means \textit{isomorphism}.

A categorical adjunction $\mathbb{L}:\mathcal{M}\leftrightarrows
\mathcal{N}:\mathbb{R}$ between two model categories is a
\textit{Quillen adjunction} if one of the following equivalent
conditions is satisfied: 1) $\mathbb{L}$ preserves cofibrations and
trivial cofibrations, 2) $\mathbb{R}$ preserves fibrations and trivial
fibrations. In that case, $\mathbb{L}$ (resp. $\mathbb{R}$) preserves
weak equivalences between cofibrant (resp. fibrant) objects.

If $P$ is a poset, let us denote by $\Delta(P)$ the \textit{order
  complex} associated with $P$. Recall that the order complex is a
simplicial complex having $P$ as underlying set and having the subsets
$\{x_0,x_1,\dots,x_n\}$ with $x_0<x_1<\dots<x_n$ as $n$-simplices
\cite{MR493916}.  Such a simplex will be denoted by
$(x_0,x_1,\dots,x_n)$.  The order complex $\Delta(P)$ can be viewed as
a poset ordered by the inclusion, and therefore as a small category.
The corresponding category will be denoted in the same way. The
opposite category $\Delta(P)^{op}$ is freely generated by the
morphisms $\de_i:(x_0,\dots,x_n) \longrightarrow
(x_0,\dots,\widehat{x_i},\dots,x_n)$ for $0\leq i\leq n$ and by the
simplicial relations $\de_i\de_j=\de_{j-1}\de_i$ for any $i<j$, where
the notation $\widehat{x_i}$ means that $x_i$ is removed.

If $\C$ is a small category, then the \textit{classifying space} of
$\C$ is denoted by $B\C$ \cite{MR0232393} \cite{MR0338129}.

\section{Reminder about the category of flows}
\label{remflow}

The category $\top$ is equipped with the unique model structure having
the \textit{weak homotopy equivalences} as weak equivalences and
having the \textit{Serre fibrations}~\footnote{that is a continuous
  map having the RLP with respect to the inclusion $\mathbf{D}^n\p
  \{0\}\subset \mathbf{D}^n\p [0,1]$ for any $n\geq 0$ where
  $\mathbf{D}^n$ is the $n$-dimensional disk.} as fibrations.

The time flow of a higher dimensional automaton is encoded in an
object called a \textit{flow} \cite{model3}. A flow $X$ consists of a
set $X^0$ called the \textit{$0$-skeleton} and whose elements
correspond to the \textit{states} (or \textit{constant execution
  paths}) of the higher dimensional automaton. For each pair of states
$(\alpha,\beta)\in X^0\p X^0$, there is a topological space
$\P_{\alpha,\beta}X$ whose elements correspond to the
\textit{(non-constant) execution paths} of the higher dimensional
automaton \textit{beginning at} $\alpha$ and \textit{ending at}
$\beta$. For $x\in \P_{\alpha,\beta}X$, let $\alpha=s(x)$ and
$\beta=t(x)$. For each triple $(\alpha,\beta,\gamma)\in X^0\p X^0\p
X^0$, there exists a continuous map $*:\P_{\alpha,\beta}X\p
\P_{\beta,\gamma}X\longrightarrow \P_{\alpha,\gamma}X$ called the
\textit{composition law} which is supposed to be associative in an
obvious sense. The topological space $\P
X=\bigsqcup_{(\alpha,\beta)\in X^0\p X^0}\P_{\alpha,\beta}X$ is called
the \textit{path space} of $X$. The category of flows is denoted by
$\dtop$. A point $\alpha$ of $X^0$ such that there are no non-constant
execution paths ending at $\alpha$ (resp. starting from $\alpha$) is
called an \textit{initial state} (resp. a \textit{final state}). A
morphism of flows $f$ from $X$ to $Y$ consists of a set map $f^0:X^0
\longrightarrow Y^0$ and a continuous map $\P f: \P X \longrightarrow
\P Y$ preserving the structure. A flow is therefore ``almost'' a small
category enriched in $\top$.

\begin{figure}
\begin{center}
\includegraphics[width=7cm]{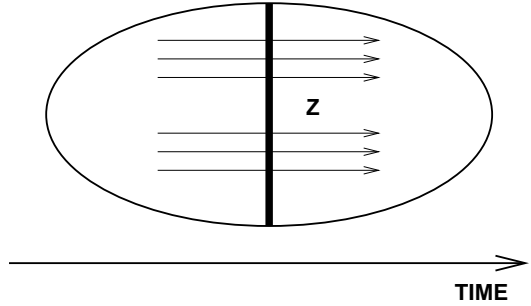}
\end{center}
\caption{Symbolic representation of
$\glob(Z)$ for some topological space $Z$} \label{exglob}
\end{figure}

An important example is the flow $\glob(Z)$ defined by the equations
\beas && \glob(Z)^0=\{\widehat{0},\widehat{1}\} \\
&& \P \glob(Z)=Z \\
&& s=\widehat{0} \\
&& t=\widehat{1} \eeas and a trivial composition law (cf.
Figure~\ref{exglob}).

The category $\dtop$ is equipped with the unique model structure
such that \cite{model3}: 
\begin{itemize}
\item The weak equivalences are the \textit{weak S-homotopy equivalences}, 
i.e. the morphisms of flows $f:X\longrightarrow Y$ such that
$f^0:X^0\longrightarrow Y^0$ is a bijection and such that $\P f:\P
X\longrightarrow \P Y$ is a weak homotopy equivalence. 
\item The fibrations are the morphisms of flows
$f:X\longrightarrow Y$ such that $\P f:\P X\longrightarrow \P Y$ is a
Serre fibration. 
\end{itemize}
This model structure is cofibrantly generated. The set of generating
cofibrations is the set $I^{gl}_+=I^{gl}\cup
\{R:\{0,1\}\longrightarrow \{0\},C:\varnothing\longrightarrow \{0\}\}$
with
\[I^{gl}=\{\glob(\mathbf{S}^{n-1})\subset \glob(\mathbf{D}^{n}), n\geq
0\}\] where $\mathbf{D}^{n}$ is the $n$-dimensional disk and
$\mathbf{S}^{n-1}$ the $(n-1)$-dimensional sphere. The set of
generating trivial cofibrations is
\[J^{gl}=\{\glob(\mathbf{D}^{n}\p\{0\})\subset
\glob(\mathbf{D}^{n}\p [0,1]), n\geq 0\}.\]

If $X$ is an object of $\cell(\dtop)$, then a presentation of the
morphism $\varnothing \longrightarrow X$ as a transfinite composition
of pushouts of morphisms of $I^{gl}_+$ is called a \textit{globular
decomposition} of $X$.

\section{Generalized T-homotopy equivalences}
\label{remindergg}

\bd A flow $X$ is {\rm loopless} if for any $\alpha\in X^0$, 
the space $\P_{\alpha,\alpha}X$ is empty. \ed

Recall that a flow is a small category without identity morphisms
enriched over a category of topological spaces.  So the preceding
definition is meaningful.

\begin{lem} \label{ordrestate} 
A flow $X$ is loopless if and only if the transitive closure of the
set $\{(\alpha,\beta)\in X^0\p X^0\hbox{ such that
}\P_{\alpha,\beta}X\neq\varnothing\}$ induces a partial ordering on
$X^0$.
\end{lem}

\bpf If $(\alpha,\beta)$ and $(\beta,\alpha)$ with $\alpha\neq \beta$
belong to the transitive closure, then there exists a finite sequence
$(x_1,\dots,x_\ell)$ of elements of $X^0$ with $x_1=\alpha$,
$x_\ell=\alpha$, $\ell>1$ and for any $m$, $\P_{x_m,x_{m+1}}X$ is
non-empty. So the space $\P_{\alpha,\alpha}X$ is non-empty because of
the existence of the composition law of $X$: contradiction.  \epf

\bd 
A {\rm full directed ball} is a flow $\vec{D}$ such that:
\begin{itemize}
\item the $0$-skeleton $\vec{D}^0$ is finite
\item $\vec{D}$ has exactly one initial state $\widehat{0}$ and one final state $\widehat{1}$  with $\widehat{0} \neq \widehat{1}$
\item each state $\alpha$ of $\vec{D}^0$ is between $\widehat{0}$ and $\widehat{1}$, that is there exists an execution path from $\widehat{0}$ to $\alpha$, and another execution path from $\alpha$ to $\widehat{1}$
\item $\vec{D}$ is loopless
\item for any $(\alpha,\beta)\in \vec{D}^0\p \vec{D}^0$, the
  topological space $\P_{\alpha,\beta}\vec{D}$ is empty if $\alpha
  \geq \beta$ and weak\-ly contractible if $\alpha<\beta$.
\end{itemize}
\ed

Let $\vec{D}$ be a full directed ball. Then by Lemma~\ref{ordrestate},
the set $\vec{D}^0$ can be viewed as a finite bounded
poset. Conversely, if $P$ is a finite bounded poset, let us consider
the \textit{flow} $F(P)$ \textit{associated with} $P$: it is of course
defined as the unique flow $F(P)$ such that $F(P)^0=P$ and
$\P_{\alpha,\beta}F(P)=\{u_{\alpha,\beta}\}$ if $\alpha<\beta$ and
$\P_{\alpha,\beta}F(P)=\varnothing$ otherwise. Then $F(P)$ is a full
directed ball and for any full directed ball $\vec{D}$, the two flows
$\vec{D}$ and $F(\vec{D}^0)$ are weakly S-homotopy equivalent.

Let $\vec{E}$ be another full directed ball. Let
$f:\vec{D}\longrightarrow\vec{E}$ be a morphism of flows preserving
the initial and final states. Then $f$ induces a morphism of posets
from $\vec{D}^0$ to $\vec{E}^0$ such that $f(\min
\vec{D}^0)=\min \vec{E}^0$ and $f(\max \vec{D}^0)=\max \vec{E}^0$. Hence 
the following definition:

\bd 
Let $\mathcal{T}$ be the class of morphisms of posets
$f:P_1\longrightarrow P_2$ such that:
\begin{enumerate}
\item The posets $P_1$ and $P_2$ are finite and bounded. 
\item The morphism of posets $f:P_1 \longrightarrow P_2$ is one-to-one; 
in particular, if $x$ and $y$ are two elements of $P_1$ with $x<y$,
then $f(x)<f(y)$.
\item One has $f(\min P_1)=\min P_2$ and  $f(\max P_1)=\max P_2$.
\end{enumerate}
Then a {\rm generalized T-homotopy equivalence} is a morphism of
$\cof(\{Q(F(f)),f\in\mathcal{T}\})$ where $Q$ is the cofibrant
replacement functor of $\dtop$.
\ed

One can choose a \textit{set} of representatives for each isomorphism
class of finite bounded posets. One obtains a \textit{set} of
morphisms $\overline{\mathcal{T}} \subset \mathcal{T}$ such that there
is the equality of classes
$\cof(\{Q(F(f)),f\in\overline{\mathcal{T}}\}) =
\cof(\{Q(F(f)),f\in\mathcal{T}\})$. By \cite{model3} Proposition~11.5,
the set of morphisms $\{Q(F(f)),f\in\overline{\mathcal{T}}\}$ permits
the small object argument.  Thus, the class of morphisms
$\cof(\{Q(F(f)),f\in\mathcal{T}\})$ contains exactly the retracts of
the morphisms of \[\cell(\{Q(F(f)),f\in\mathcal{T}\})\] by
\cite{MR99h:55031} Corollary~2.1.15.

The inclusion of posets $\{\widehat{0} < \widehat{1}\} \subset
\{\widehat{0} < 2 < \widehat{1}\}$ corresponds to the case of
Figure~\ref{ex1}.

\begin{figure}
\begin{center}
\includegraphics[width=9cm]{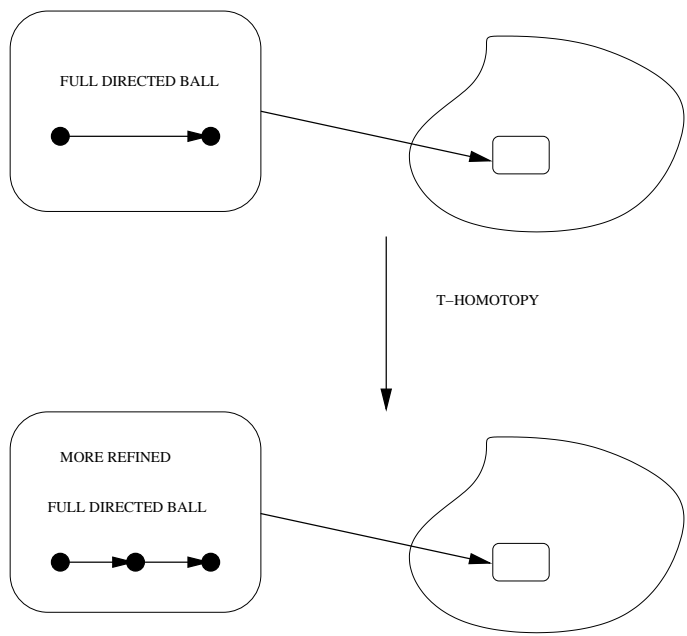}
\end{center}
\caption{Replacement of a full directed ball by a more refined one}
\label{ex2}
\end{figure}

A T-homotopy consists in locally replacing in a flow a full directed
ball by a more refined one (cf. Figure~\ref{ex2}), and in iterating
the process transfinitely.

\begin{figure}
\[
\xymatrix{\widehat{0} \ar@{->}[rrrr]^-{U} &&&& \widehat{1} \\
\widehat{0} \ar@{->}[rr]^-{U'} && A \ar@{->}[rr]^-{U''} && \widehat{1}}
\]
\caption{The simplest example of refinement of observation}
\label{ex1}
\end{figure}

\section{Globular complex}
\label{remindergl}

The technical reference is \cite{model2}. A \textit{globular complex}
is a topological space together with a structure describing the
sequential process of attaching \textit{globular cells}. A general
globular complex may require an arbitrary long transfinite
construction. We restrict our attention in this paper to globular
complexes whose globular cells are morphisms of the form
$\glob^{top}(\mathbf{S}^{n-1}) \longrightarrow
\glob^{top}(\mathbf{D}^{n})$.

\bd 
A {\rm multipointed topological space} $(X,X^0)$ is a pair of
topological spaces such that $X^0$ is a discrete subspace of $X$.  A
morphism of multipointed topological spaces $f:(X,X^0)\longrightarrow
(Y,Y^0)$ is a continuous map $f:X\longrightarrow Y$ such that
$f(X^0)\subset Y^0$. The corresponding category is denoted by
$\top^m$. The set $X^0$ is called the {\rm $0$-skeleton} of $(X,X^0)$.
The space $X$ is called the {\rm underlying topological space} of
$(X,X^0)$. 
\ed

The category of multipointed spaces is cocomplete.

\bd 
Let $Z$ be a topological space. The {\rm globe of $Z$}, which is
denoted by $\glob^{top}(Z)$, is the multipointed space
\[(|\glob^{top}(Z)|,\{\widehat{0},\widehat{1}\})\] 
where the topological space $|\glob^{top}(Z)|$ is the quotient of
$\{\widehat{0},\widehat{1}\}\sqcup (Z\p[0,1])$ by the relations
$(z,0)=(z',0)=\widehat{0}$ and $(z,1)=(z',1)=\widehat{1}$ for any
$z,z'\in Z$. In particular, $\glob^{top}(\varnothing)$ is the
multipointed space
$(\{\widehat{0},\widehat{1}\},\{\widehat{0},\widehat{1}\})$.
\ed

If $Z$ is not empty, then the space $|\glob^{top}(Z)|$ is the
unpointed suspension of $Z$. If $Z$ is the empty space, then the space
$|\glob^{top}(Z)|$ is the discrete two-point space.

\begin{nota} 
Let  $Z$ be a singleton. The globe of $Z$ is denoted by
$\vI^{top}$. 
\end{nota}

\bd 
Let $I^{gl,top}:=\{\glob^{top}(\mathbf{S}^{n-1})\longrightarrow
\glob^{top}(\mathbf{D}^{n}),n\geq 0\}$. 
A {\rm relative globular precomplex} is a relative $I^{gl,top}$-cell
complex in the category of multipointed topological spaces. 
\ed

\bd 
A {\rm globular precomplex} is a $\lambda$-sequence of multipointed
topological spaces $X:\lambda\longrightarrow \top^m$ such that $X$ is
a relative globular precomplex and such that $X_0=(X^0,X^0)$ with
$X^0$ a discrete space. This $\lambda$-sequence is characterized by a
presentation ordinal $\lambda$, and for any $\beta<\lambda$, an
integer $n_\beta\geq 0$ and an attaching map $\phi_\beta :
\glob^{top}(\mathbf{S}^{n_\beta-1}) \longrightarrow X_\beta$. The
family $(n_\beta,\phi_\beta)_{\beta<\lambda}$ is called the {\rm
globular decomposition} of $X$.
\ed

Let $X$ be a globular precomplex. The $0$-skeleton of $\liminj X$ is
equal to $X^0$.

\bd A morphim of globular precomplexes $f:X\longrightarrow Y$ is a
morphism of multipointed spaces still denoted by $f$ from $\liminj X$
to $\liminj Y$. \ed

\begin{nota} 
If $X$ is a globular precomplex, then the underlying topological space
of the multipointed space $\liminj X$ is denoted by $|X|$ and the
$0$-skeleton of the multipointed space $\liminj X$ is denoted by
$X^0$. 
\end{nota}

\bd 
Let $X$ be a globular precomplex. The space $|X|$ is called the {\rm
underlying topological space} of $X$. The set $X^0$ is called the {\rm
$0$-skeleton} of $X$. 
\ed

\bd 
Let $X$ be a globular precomplex. A morphism of globular precomplexes
$\gamma:\vI^{top}\longrightarrow X$ is a {\rm non-constant execution
path} of $X$ if there exists $t_0=0<t_1<\dots<t_{n}=1$ such that:
\begin{enumerate}
\item  $\gamma(t_i)\in X^0$ for any $0 \leq i \leq n$, 
\item $\gamma(]t_i,t_{i+1}[)\subset \glob^{top}(\mathbf{D}^{n_{\beta_i}} 
\backslash \mathbf{S}^{n_{\beta_i}-1})$ for some $(n_{\beta_i},\phi_{\beta_i})$ 
of the globular decomposition of $X$, 
\item  for $0\leq i<n$, there exists $z^i_\gamma\in \mathbf{D}^{n_{\beta_i}}\backslash \mathbf{S}^{n_{\beta_i}-1}$ and a strictly increasing continuous map
$\psi^i_\gamma:[t_i,t_{i+1}]\longrightarrow [0,1]$ such that
$\psi^i_\gamma(t_i)=0$ and $\psi^i_\gamma(t_{i+1})=1$ and for any
$t\in [t_i,t_{i+1}]$, $\gamma(t)=(z^i_\gamma,\psi^i_\gamma(t))$.
\end{enumerate}
In particular, the restriction $\gamma\!\restriction_{]t_i,t_{i+1}[}$
of $\gamma$ to $]t_i,t_{i+1}[$ is one-to-one. The set of non-constant
execution paths of $X$ is denoted by ${\P}^{top}(X)$. 
\ed

\bd 
A morphism of globular precomplexes $f:X\longrightarrow Y$ is {\rm
non-decreasing} if the canonical set map
$\top([0,1],|X|)\longrightarrow \top([0,1],|Y|)$ induced by
composition by $f$ yields a set map ${\P}^{top}(X)\longrightarrow
{\P}^{top}(Y)$. In other terms, one has the commutative diagram of
sets
\[\xymatrix{
{\P}^{top}(X)\fr{}\fd{\subset}& {\P}^{top}(Y)\fd{\subset}\\
\top([0,1],|X|) \fr{} &\top([0,1],|Y|).}
\]
\ed

\bd 
A {\rm globular complex } (resp. a {\rm relative globular complex})
$X$ is a globular precomplex (resp. a relative globular precomplex)
such that the attaching maps $\phi_\beta$ are non-decreasing. A
morphism of globular complexes is a morphism of globular precomplexes
which is non-decreasing. The category of globular complexes together
with the morphisms of globular complexes as defined above is denoted
by $\gltop$. 
\ed

\bd Let $X$ be a globular complex. A point $\alpha$ of $X^0$ such that
there are no non-constant execution paths ending at $\alpha$ (resp.
starting from $\alpha$) is called {\rm initial state} (resp.  {\rm
  final state}). More generally, a point of $X^0$ will be sometime
called {\rm a state} as well.  \ed

\bth (\cite{model2} Theorem~III.3.1) 
There exists a unique functor $cat:\gltop\longrightarrow\dtop$ such
that
\begin{enumerate}
\item if $X=X^0$ is a discrete globular complex, then $cat(X)$ is
the achronal flow $X^0$ (``achronal'' meaning with an empty path space)
\item if $Z=\mathbf{S}^{n-1}$ or $Z=\mathbf{D}^{n}$ for some integer $n\geq 0$, 
then $cat(\glob^{top}(Z))=\glob(Z)$, 
\item for any globular complex $X$ with globular decomposition
$(n_\beta,\phi_\beta)_{\beta<\lambda}$, for any limit
ordinal $\beta\leq\lambda$, the canonical morphism of flows
\[\liminj_{\alpha<\beta} cat(X_\alpha)\longrightarrow
cat(X_\beta)\] is an isomorphism of flows, 
\item for any globular complex $X$ with globular decomposition
$(n_\beta,\phi_\beta)_{\beta<\lambda}$, for any
$\beta<\lambda$, one has the pushout of flows
\[\xymatrix{\glob(\mathbf{S}^{n_\beta-1})\fr{cat(\phi_\beta)}\fd{}& cat(X_\beta)\fd{}\\
\glob(\mathbf{D}^{n_\beta})\fr{} & cat(X_{\beta+1}).\cocartesien}\]
\end{enumerate}
\eth

\section{The underlying homotopy type of a flow}
\label{defunder}

\bth \label{remonte}
The functor $cat$ induces a functor, still denoted by $cat$ from
$\gltop$ to $\hda$. For any flow $X$ of $\hda$, there exists a
globular complex $Y$ such that $cat(Y)=X$. It is constructed by using
the globular decomposition of $X$. If two globular complexes $Y_1$ and
$Y_2$ satisfy $cat(Y_1)=cat(Y_2)=X$, then the two topological spaces
$|Y_1|$ and $|Y_2|$ are homotopy equivalent.
\eth

\bpf 
The construction of $Y$ is made in the proof of \cite{model2}
Theorem~V.4.1.  If two globular complexes $Y_1$ and $Y_2$ satisfy
$cat(Y_1)=cat(Y_2)=X$, then they are S-homotopy equivalent by
\cite{model2} Theorem~IV.4.9. And the S-homotopy equivalence between the 
globular complexes $Y_1$ and $Y_2$ yields an homotopy equivalence
between the underlying topological spaces $|Y_1|$ and $|Y_2|$ by
\cite{model2} Proposition~VII.2.2.
\epf

The recipe to obtain the underlying homotopy type of a flow $X$ is as
follows \cite{model2}:
\begin{enumerate}
\item Take a flow $X$.
\item Take its cofibrant replacement $Q(X)\in\hda$.
\item By Theorem~\ref{remonte}, there exists a globular complex 
$X^{top}$ such that $cat(X^{top})=Q(X)$.
\item The cofibrant topological space $|X^{top}|$ is unique up to 
homotopy and is the underlying homotopy type $|X|$ of $X$.
\end{enumerate}
This yields a well defined functor $|-|:\dtop\longrightarrow
\ho(\top)$ from the category of flows to the homotopy category of
topological spaces (\cite{model2} Part~VII.2).

Roughly speaking, the underlying homotopy type of a flow $X$ consists
in factoring the morphism of flows $X^0\longrightarrow Q(X)$ as a
transfinite composition of pushouts of elements of $I^{gl}$; and then
replacing this transfinite composition by a transfinite composition of
pushouts of the continuous maps $\{|\glob^{top}(\mathbf{S}^{n-1})|
\subset |\glob^{top}(\mathbf{D}^{n})|, n\geq 0\}$; and then
calculating this transfinite composition in $\top$: the result is a
cofibrant topological space which is unique up to homotopy.

\section{A useful Reedy category and related lemmas}
\label{reedyuse}

Let $P$ be a finite bounded poset with bottom element $\widehat{0}$
and with top element $\widehat{1}$.  Let us denote by
$\Delta^{ext}(P)$ the full subcategory of $\Delta(P)$ consisting of
the simplices $(\alpha_0,\dots,\alpha_p)$ such that $\widehat{0} =
\alpha_0$ and $\widehat{1} = \alpha_p$.
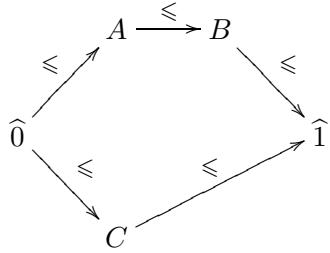
\begin{figure}
\[
\xymatrix{
& A \fr{\leq} & B \ar@{->}[rd]^{\leq} & \\
\widehat{0} \ar@{->}[rd]^{\leq} \ar@{->}[ru]^{\leq} & & & \widehat{1}\\
&  C \ar@{->}[rru]^{\leq} & & }
\] 
\caption{Example of finite bounded poset}
\label{poset}
\end{figure}
If $P=\{\widehat{0}<A<B<\widehat{1}, \widehat{0}<C<\widehat{1}\}$ is
the poset of Figure~\ref{poset}, then the small category
$\Delta^{ext}(P)^{op}$ looks as follows:
\[
\xymatrix{
& (\widehat{0},A,B,\widehat{1})\ar@{->}[d]\ar@{->}[dr] & & (\widehat{0},C,\widehat{1})\ar@/^30pt/[ldd] & \\
& (\widehat{0},A,\widehat{1})\ar@{->}[dr] & (\widehat{0},B,\widehat{1})\ar@{->}[d] & & \\
& & (\widehat{0},\widehat{1})   & & 
}
\]
The simplex $(\widehat{0},\widehat{1})$ is always a terminal object of
$\Delta^{ext}(P)^{op}$.

\begin{nota} \label{defell} 
Let $X$ be a loopless flow such that $(X^0,\leq)$ is locally finite. 
Let $(\alpha,\beta)$ be a $1$-simplex of $\Delta(X^0)$. We denote by
$\ell(\alpha,\beta)$ the maximum of the set of integers
\[\left\{p\geq 1, \exists (\alpha_0,\dots,\alpha_p)\hbox{ $p$-simplex of }\Delta(X^0)\hbox{ s.t. }
(\alpha_0,\alpha_p)= (\alpha,\beta)\right\}\] One always has $1\leq
\ell(\alpha,\beta)\leq \card(]\alpha,\beta])$.
\end{nota}

\begin{lem}\label{triinv} 
Let $X$ be a loopless flow such that $(X^0,\leq)$ is locally finite.
Let $(\alpha,\beta,\gamma)$ be a $2$-simplex of $\Delta(X^0)$. Then
one has
\[\ell(\alpha,\beta)+\ell(\beta,\gamma)\leq \ell(\alpha,\gamma).\] 
\end{lem}

\bpf 
Let $\alpha = \alpha_0 < \dots < \alpha_{\ell(\alpha,\beta)} =
\beta$. Let $\beta = \beta_0 < \dots < \beta_{\ell(\beta,\gamma)} = 
\gamma$. Then 
\[(\alpha_0,\dots,\alpha_{\ell(\alpha,\beta)},\beta_1,\dots,\beta_{\ell(\beta,\gamma)}) \]
is a simplex of $\Delta(X^0)$ with $\alpha = \alpha_0$ and $\beta_{\ell(\beta,\gamma)} = 
\gamma$. So $\ell(\alpha,\beta)+\ell(\beta,\gamma)\leq \ell(\alpha,\gamma)$. 
\epf

\bp \label{ReedyTwo}
Let $P$ be a finite bounded poset.  Let
\[
d(\alpha_0,\dots,\alpha_p)=\ell(\alpha_0,\alpha_1)^2+\dots + 
\ell(\alpha_{p-1},\alpha_p)^2\]
where $\ell$ is the function of Notation~\ref{defell}. Then $d$ yields
a functor $\Delta^{ext}(P)^{op}\longrightarrow \N$ making $\Delta^{ext}(P)^{op}$
a direct category, that is a Reedy category with $\Delta^{ext}(P)^{op}_+: =
\Delta^{ext}(P)^{op}$ and $\Delta^{ext}(P)^{op}_- = \varnothing$.
\ep

\bpf 
Let $\de_i:(\alpha_0,\dots,\alpha_p)\longrightarrow
(\alpha_0,\dots,\widehat{\alpha_i},\dots,\alpha_p)$ be a morphism of
$\Delta^{ext}(P)^{op}$ with $p\geq 2$ and $0< i<p$ .  Then
\beas
&& d(\alpha_0,\dots,\alpha_p) = \ell(\alpha_0,\alpha_1)^2+\dots+\ell(\alpha_{p-1},\alpha_p)^2\\
&& d(\alpha_0,\dots,\widehat{\alpha_i},\dots,\alpha_p) =\ell(\alpha_0,\alpha_1)^2+\dots+\ell(\alpha_{i-1},\alpha_{i+1})^2+\dots+\ell(\alpha_{p-1},\alpha_p)^2.
\eeas
So one obtains 
\[ d(\alpha_0,\dots,\alpha_p)-d(\alpha_0,\dots,\widehat{\alpha_i},\dots,\alpha_p)
=\ell(\alpha_{i-1},\alpha_{i})^2+\ell(\alpha_{i},\alpha_{i+1})^2-\ell(\alpha_{i-1},\alpha_{i+1})^2.\]
By Lemma~\ref{triinv}, one has 
\[(\ell(\alpha_{i-1},\alpha_{i})+\ell(\alpha_{i},\alpha_{i+1}))^2\leq
\ell(\alpha_{i-1},\alpha_{i+1})^2.\] Therefore, one obtains 
\[\ell(\alpha_{i-1},\alpha_{i})^2+\ell(\alpha_{i},\alpha_{i+1})^2<\ell(\alpha_{i-1},\alpha_{i+1})^2
\]
since $2\ell(\alpha_{i-1},\alpha_{i})\ell(\alpha_{i},\alpha_{i+1})\geq
2$. Thus, every morphism of $\Delta^{ext}(P)^{op}$ raises the degree.
\epf

\begin{cor} \label{calculbienholim}
Let $P$ be a finite bounded poset. Then the colimit functor
\[\liminj:\top^{\Delta^{ext}(P)^{op}\backslash
\{(\widehat{0},\widehat{1})\}}\longrightarrow \top\] 
is a left Quillen functor if the category of diagrams
$\top^{\Delta^{ext}(P)^{op}\backslash
\{(\widehat{0},\widehat{1})\}}$ is equipped with the Reedy model
structure.
\end{cor}

Indeed, the fact that the colimit functor is a left Quillen functor
will be actually applied for $\Delta^{ext}(P)^{op}\backslash
\{(\widehat{0},\widehat{1})\}$. Recall that the pair
$(\widehat{0},\widehat{1})$ is a terminal object of
$\Delta^{ext}(P)^{op}$. Therefore, it is not particularly interesting
to calculate the colimit of a diagram of spaces over the whole
category $\Delta^{ext}(P)^{op}$. Note also that there is an
isomorphism of small categories
\[\boxed{\Delta^{ext}(P)^{op}\backslash
\{(\widehat{0},\widehat{1})\} \iso \de(\Delta^{ext}(\vec{D}^0)^{op}_+\!\downarrow\! (\widehat{0},\widehat{1}))}.\]

\bpf The Reedy structure on $\Delta^{ext}(P)^{op}\backslash
\{(\widehat{0},\widehat{1})\}$ provides a model structure on the
category $\top^{\Delta^{ext}(P)^{op}\backslash
  \{(\widehat{0},\widehat{1})\}}$ of diagrams of topological spaces
over the small category \[\Delta^{ext}(P)^{op}\backslash
\{(\widehat{0},\widehat{1})\}\] such that a morphism of diagrams
$f:D\longrightarrow E$ is
\begin{enumerate}
\item a weak equivalence if and only if for every object
  ${\underline{\alpha}}$ of $\Delta^{ext}(P)^{op}\backslash
  \{(\widehat{0},\widehat{1})\}$, the morphism
  $D_{\underline{\alpha}}\longrightarrow E_{\underline{\alpha}}$ is a
  weak equivalence of $\top$ (we will use the term objectwise weak
  equivalence to describe this situation)
\item a cofibration if and only if for every object
  ${\underline{\alpha}}$ of $\Delta^{ext}(P)^{op}\backslash
  \{(\widehat{0},\widehat{1})\}$, the morphism
  $D_{\underline{\alpha}}\sqcup_{L_{\underline{\alpha}} D}
  L_{\underline{\alpha}}E \longrightarrow E_{\underline{\alpha}}$ is a
  cofibration of $\top$
\item a fibration if and only if for every object
  ${\underline{\alpha}}$ of $\Delta^{ext}(P)^{op}\backslash
  \{(\widehat{0},\widehat{1})\}$, the morphism
  $D_{\underline{\alpha}}\longrightarrow
  E_{\underline{\alpha}}\p_{M_{\underline{\alpha}} E}
  M_{\underline{\alpha}}D$ is a fibration of $\top$.
\end{enumerate}
For every object $\underline{\alpha}$ of
$\Delta^{ext}(P)^{op}\backslash \{(\widehat{0},\widehat{1})\}$, the
matching category
$\de(\underline{\alpha}\!\downarrow\!\Delta^{ext}(P)^{op}_-)$ is
empty.  So for every object $A$ of the diagram category
$\top^{\Delta^{ext}(P)^{op}}$ and every object $\underline{\alpha}$ of
the small category $\Delta^{ext}(P)^{op}\backslash
\{(\widehat{0},\widehat{1})\}$, there is an isomorphism
$M_{\underline{\alpha}} A\iso \mathbf{1}$.  So a Reedy fibration is an
objectwise fibration.  Therefore, the diagonal functor $\diag$ of the
adjunction $\liminj:\top^{\Delta^{ext}(P)^{op} \backslash
  \{(\widehat{0},\widehat{1})\}} \leftrightarrows \top:\diag$ is a
right Quillen functor.  \epf

\bp \label{constructionhat}
Let $\vec{D}$ be a full directed ball. There exists one and only one
functor
\[\mathcal{F}{\vec{D}}:\Delta^{ext}(\vec{D}^0)^{op}\longrightarrow \top\]
satisfying the following conditions:
\begin{enumerate}
\item 
$\mathcal{F}{\vec{D}}_{(\alpha_0,\dots,\alpha_p)} = \P_{\alpha_0,\alpha_1}\vec{D}
\p\dots \p \P_{\alpha_{p-1},\alpha_p}\vec{D}$ (recall that necessarily,
one has the equalities $\alpha_0=\widehat{0}$ and $\alpha_p=\widehat{1}$ by definition 
of the small category $\Delta^{ext}(\vec{D}^0)^{op}$)
\item 
the unique morphism
$\de_i:\mathcal{F}{\vec{D}}_{(\alpha_0,\dots,\alpha_p)}\longrightarrow
\mathcal{F}{\vec{D}}_{(\alpha_0,\dots,\widehat{\alpha_i},\dots,\alpha_p)}$ for
$0<i<p$ is induced by the composition law
$\P_{\alpha_{i-1},\alpha_i}\vec{D}\p
\P_{\alpha_{i},\alpha_{i+1}}\vec{D}\longrightarrow
\P_{\alpha_{i-1},\alpha_{i+1}}\vec{D}$. 
\end{enumerate}
\ep

\bpf The uniqueness on objects is exactly the first assertion. 
The uniqueness on morphisms comes from the fact that any morphism of
$\Delta^{ext}(\vec{D}^0)^{op}$ is a composite of $\de_i$. We have to
prove the existence. The diagram of topological spaces
\[
\xymatrix{
\mathcal{F}{\vec{D}}_{(\alpha_0,\dots,\alpha_p)} \fr{\de_i}\fd{\de_j} & \mathcal{F}{\vec{D}}_{(\alpha_0,\dots,\widehat{\alpha_i},\dots,\alpha_p)}\fd{\de_{j-1}}\\
\mathcal{F}{\vec{D}}_{(\alpha_0,\dots,\widehat{\alpha_j},\dots,\alpha_p)}\fr{\de_i}& \mathcal{F}{\vec{D}}_{(\alpha_0,\dots,\widehat{\alpha_i},\dots,\widehat{\alpha_j},\dots,\alpha_p)}}
\]
is commutative for any $0<i<j<p$ and any $p\geq 2$. Indeed, if $i<j-1$, 
then one has
\[\boxed{\de_i\de_j(\gamma_1,\dots,\gamma_{p})=\de_{j-1}\de_i(\gamma_1,\dots,\gamma_{p})=
(\gamma_1,\dots,\gamma_{i}\gamma_{i+1},\dots,\gamma_{j}\gamma_{j+1},\dots,\gamma_{p})}\]
and if $i=j-1$, then one has 
\[\boxed{\de_i\de_j(\gamma_1,\dots,\gamma_{p})=\de_{j-1}\de_i(\gamma_1,\dots,\gamma_{p})=
(\gamma_1,\dots,\gamma_{j-1}\gamma_{j}\gamma_{j+1},\dots,\gamma_{p})}\]
because of the associativity of the composition law of $X$.  In other
terms, the $\de_i$ maps satisfy the simplicial identities. Hence the
result. 
\epf

Take again the poset $P$ of Figure~\ref{poset}, and the corresponding
full directed ball $\vec{D}=F(P)$. The diagram $\mathcal{F}{\vec{D}}$
looks as follows: {\small \[
\xymatrix{
& \P_{\widehat{0},A}\vec{D}\p \P_{A,B}\vec{D}\p \P_{B,\widehat{1}}\vec{D}\ar@{->}[d]\ar@{->}[dr] & & \P_{\widehat{0},C}\vec{D}\p \P_{C,\widehat{1}}\vec{D}\ar@/^30pt/[ldd] & \\
& \P_{\widehat{0},A}\vec{D}\p \P_{A,\widehat{1}}\vec{D}\ar@{->}[dr] & \P_{\widehat{0},B}\vec{D}\p \P_{B,\widehat{1}}\vec{D}\ar@{->}[d] & & \\
& & \P_{\widehat{0},\widehat{1}}\vec{D}   & & 
}
\]}

\bd Let $X$ be a flow. Let $A$ be a subset of $X^0$. Then the
\textit{restriction} $X\!\restriction_A$ of $X$ over $A$ is the unique
flow such that $(X\!\restriction_A)^0=A$, such that
$\P_{\alpha,\beta}(X\!\restriction_A)=\P_{\alpha,\beta}X$ for
$(\alpha,\beta)\in A\p A$ and such that the inclusions $A\subset X^0$
and $\P (X\!\restriction_A)\subset \P X$ induce a morphism of flows
$X\!\restriction_A\longrightarrow X$.  \ed

\bp Let $\vec{D}$ be a full directed ball. Let $(\alpha,\beta)$ be a
simplex of $\Delta(\vec{D}^0)$. Then
$\vec{D}\!\restriction_{[\alpha,\beta]}$ is a full directed ball with
initial state $\alpha$ and with final state $\beta$.  \ep

\bpf Obvious. \epf

\bp 
Let $\vec{D}$ and $\vec{D'}$ be two full directed balls. Then the flow
$\vec{D}*\vec{D'}$ obtained by identifying the final state
$\widehat{1}$ of $\vec{D}$ with the initial state $\widehat{0}$ of
$\vec{D'}$ is a full directed ball.
\ep

\bpf The condition which is less easy to verify than the other ones
is: for any $(\alpha,\beta)\in (\vec{D}*\vec{D'})^0 \p
(\vec{D}*\vec{D'})^0$, the topological space $\P_{\alpha,\beta}
(\vec{D}*\vec{D'})$ is weakly contractible if $\alpha<\beta$. Let $m$
be the point of $\vec{D}*\vec{D'}$ corresponding to the final state of
$\vec{D}$ and the initial state of $\vec{D'}$. If $\alpha<\beta\leq
m$, then one has the isomorphism of spaces $\P_{\alpha,\beta}
(\vec{D}*\vec{D'}) \iso \P_{\alpha,\beta} \vec{D}$. If $m\leq
\alpha<\beta$, then one has the isomorphism of spaces
$\P_{\alpha,\beta} (\vec{D}*\vec{D'}) \iso \P_{\alpha,\beta}
\vec{D'}$. At last, if $\alpha<m<\beta$, then one has the isomorphism
of spaces $\P_{\alpha,\beta} (\vec{D}*\vec{D'}) \iso \P_{\alpha,m}
\vec{D} \p \P_{m,\beta} \vec{D'}$. So in each case, the space
$\P_{\alpha,\beta} (\vec{D}*\vec{D'})$ is weakly contractible.  \epf

\bp Let $\vec{D}$ be a full directed ball.  There exists one and only one functor 
\[\mathcal{G}{\vec{D}}:\Delta^{ext}(\vec{D}^0)^{op}\longrightarrow \dtop\]
satisfying the following conditions:
\begin{enumerate}
\item for any object $(\alpha_0,\dots,\alpha_p)$ of $\Delta^{ext}(\vec{D}^0)^{op}$, let 
\[\mathcal{G}{\vec{D}}(\alpha_0,\dots,\alpha_p)=\vec{D}\!\restriction_{[\alpha_0,\alpha_1]} * \dots *
\vec{D}\!\restriction_{[\alpha_{p-1},\alpha_p]}\] 
\item the unique morphism
$\mathcal{G}{\vec{D}}(\alpha_0,\dots,\alpha_p)\longrightarrow
\mathcal{G}{\vec{D}}(\alpha_0,\dots,\widehat{\alpha_i},\dots,\alpha_p)$ for
$0<i<p$ is induced by the composition law
$\vec{D}\!\restriction_{[\alpha_{i-1},\alpha_i]}*
\vec{D}\!\restriction_{[\alpha_{i},\alpha_{i+1}]}\longrightarrow
\vec{D}\!\restriction_{[\alpha_{i-1},\alpha_{i+1}]}$. 
\end{enumerate}
\ep

Notice that $\vec{D}\!\restriction_{[\widehat{0},\widehat{1}]}=\vec{D}$.

\bpf This comes from the associativity of the composition law of a flow. 
\epf

\bp 
Let $\vec{D}$ be a full directed ball.  Let
$(\alpha_0,\dots,\alpha_p)\in \Delta^{ext}(\vec{D}^0)^{op}$ be a
simplex. Then there exists a unique morphism of flows
{\small \[u_{(\alpha_0,\dots,\alpha_p)}:\glob(\P_{\alpha_0,\alpha_1}\vec{D}
\p\dots \p \P_{\alpha_{p-1},\alpha_p}\vec{D})\longrightarrow 
\glob(\P_{\alpha_0,\alpha_1}\vec{D}) * 
\dots * \glob(\P_{\alpha_{p-1},\alpha_p}\vec{D})\]} 
such that $u_{(\alpha_0,\dots,\alpha_p)}(x_1,\dots,x_p)=x_1*\dots
*x_p$. With $(\alpha_0,\dots,\alpha_p)$ running over the set of
simplices of $\Delta^{ext}(\vec{D}^0)^{op}$, one obtains a morphism of
diagrams of flows \[\glob(\mathcal{F}{\vec{D}})\longrightarrow
\mathcal{G}{\vec{D}}.\] 
\ep

\bpf Obvious. 
\epf

\bp \label{pushmax}
Let $\vec{D}$ be a full directed ball.  Then one has the pushout diagram of flows: 
\[
\xymatrix{
\glob(L_{(\widehat{0},\widehat{1})}\mathcal{F}{\vec{D}})\fr{}\fd{} & L_{(\widehat{0},\widehat{1})}\mathcal{G}{\vec{D}}\fd{}\\
\glob(\P_{\widehat{0},\widehat{1}}\vec{D}) \fr{} & \vec{D}. \cocartesien}
\]
\ep

This statement remains true when the $1$-simplex
$(\widehat{0},\widehat{1})$ is replaced by another $1$-simplex
$(\alpha,\beta)$ of $\Delta(\vec{D}^0)$. This statement above becomes
false in general when the $1$-simplex $(\widehat{0},\widehat{1})$ is
replaced by a $p$-simplex of $\Delta(\vec{D}^0)$ with $p\geq 2$.

Let us illustrate this proposition in the case of $\vec{D}^0 =
\{\widehat{0} < A < \widehat{1}\}$. One then has: 
\begin{enumerate}
\item $L_{(\widehat{0},\widehat{1})}\mathcal{F}{\vec{D}} = \P_{\widehat{0},A}\vec{D} \p \P_{A,\widehat{1}}\vec{D}$; 
\item $L_{(\widehat{0},\widehat{1})}\mathcal{G}{\vec{D}} = \vec{D}\!\restriction_{[\widehat{0},A]}*  \vec{D}\!\restriction_{[A,\widehat{1}]} = \glob(\P_{\widehat{0},A}\vec{D}) * \glob(\P_{A,\widehat{1}}\vec{D})$; the last equality is due to the fact that $]\widehat{0},A[ = ]A,\widehat{1}[ = \varnothing$.
\item The pushout above is equivalent to the following pushout: 
\[
\xymatrix{
\glob(\P_{\widehat{0},A}\vec{D} \p \P_{A,\widehat{1}}\vec{D})\fr{}\fd{} & \glob(\P_{\widehat{0},A}\vec{D}) * \glob(\P_{A,\widehat{1}}\vec{D})\fd{}\\
\glob(\P_{\widehat{0},\widehat{1}}\vec{D}) \fr{} & \vec{D}. \cocartesien}
\]
\end{enumerate}

\bpf One already has the commutative diagram 
\[
\xymatrix{
\glob(L_{(\widehat{0},\widehat{1})}\mathcal{F}{\vec{D}})\fr{}\fd{} & L_{(\widehat{0},\widehat{1})}\mathcal{G}{\vec{D}}\fd{}\\
\glob(\P_{\widehat{0},\widehat{1}}\vec{D}) \fr{} & \vec{D}.}
\]
Therefore, one only has to check that $\vec{D}$ satisfies the same
universal property as the pushout.

Consider a commutative diagram of flows of the form: 
\[
\xymatrix{
\glob(L_{(\widehat{0},\widehat{1})}\mathcal{F}{\vec{D}})\fr{}\fd{} & L_{(\widehat{0},\widehat{1})}\mathcal{G}{\vec{D}}\fd{}\\
\glob(\P_{\widehat{0},\widehat{1}}\vec{D}) \fr{} & X.}
\]
The morphism of flows $\glob(\P_{\widehat{0},\widehat{1}}\vec{D})
\longrightarrow X$ induces a continuous map
$\P_{\widehat{0},\widehat{1}}\vec{D} \longrightarrow \P X$. The
morphism of flows $L_{(\widehat{0},\widehat{1})}\mathcal{G}{\vec{D}}
\longrightarrow X$ induces a continuous map
$\P_{\alpha,\beta}\vec{D}\longrightarrow \P X$ for any $1$-simplex
$(\alpha,\beta)$ of $\Delta(\vec{D}^0)$ with $(\alpha,\beta)\neq
(\widehat{0},\widehat{1})$. The existence of the morphism of flows
$L_{(\widehat{0},\widehat{1})}\mathcal{G}{\vec{D}} \longrightarrow X$
ensures the compatibility of the continuous maps
$\P_{\alpha,\beta}\vec{D}\longrightarrow \P X$ for $(\alpha,\beta) \in
\Delta(\vec{D}^0)$ with the composition of execution paths involving a
triple $(\alpha,\beta,\gamma)$ such that $(\alpha,\gamma)
\neq(\widehat{0},\widehat{1})$. And the commutativity of the diagram
with $X$ ensures the compatibility of the continuous maps
$\P_{\alpha,\beta}\vec{D}\longrightarrow \P X$ for $(\alpha,\beta) \in
\Delta(\vec{D}^0)$ with the composition of execution paths involving a
triple $(\alpha,\beta,\gamma)$ such that $(\alpha,\gamma)
=(\widehat{0},\widehat{1})$. Hence the existence and uniqueness of the
morphism $\vec{D}\longrightarrow X$. \epf

\bth (\cite{MR1712872} Theorem~1 p.~213) \label{souslim}
Let $L:J'\longrightarrow J$ be a final functor between small
categories, i.e.  such that for any $k\in J$, the comma category
$(k\!\downarrow\!L)$ is non-empty and connected. Let
$F:J\longrightarrow \C$ be a functor from $J$ to a cocomplete category
$\C$. Then $L$ induces a canonical morphism $\liminj F\circ
L\longrightarrow \liminj F$ which is an isomorphism.
\eth

\begin{nota} Let $X$ be a loopless flow. 
Let $\underline{\alpha}=(\alpha_0,\dots,\alpha_p)$ be a simplex of the
order complex $\Delta(X^0)$ of the poset $X^0$. Let
$\alpha<\alpha_0$. Then the notation $\alpha.\underline{\alpha}$
represents the simplex $(\alpha,\alpha_0,\dots,\alpha_p)$ of
$\Delta(X^0)$.
\end{nota}

\bth \label{simplification}
Let $\vec{D}$ be a full directed ball. Let $\underline{\alpha} =
(\alpha_0,\dots,\alpha_p)$ be a simplex of
$\Delta^{ext}(\vec{D}^0)^{op}$. Let
$i_{(\alpha_0,\dots,\alpha_p)}:L_{(\alpha_0,\dots,\alpha_p)}\mathcal{F}\vec{D}
\longrightarrow \mathcal{F}\vec{D}_{(\alpha_0,\dots,\alpha_p)}$.  Then
one has \[i_{(\alpha_0,\dots,\alpha_p)} =
i_{(\alpha_0,\alpha_1)}\square \dots \square
i_{(\alpha_{p-1},\alpha_p)}\] where $\square$ is the pushout product
(cf. Notation~\ref{notapushcal}).
\eth

\bpf Let $\underline{\alpha} = (\alpha_0,\dots,\alpha_p)$ be a fixed
object of $\Delta^{ext}(\vec{D}^0)^{op}$. The latching category
\[\de(\Delta^{ext}(\vec{D}^0)^{op}_+\!\downarrow\!
\underline{\alpha})\] is the full subcategory of
$\Delta^{ext}(\vec{D}^0)^{op}$ consisting of the simplices
${\underline{\beta} = (\beta_0,\dots,\beta_q)}$ such that there is a
strict inclusion \[\{\alpha_0,\dots,\alpha_p\} \subsetneqq
\{\beta_0,\dots,\beta_q\},\] that is $\{\alpha_0,\dots,\alpha_p\}
\subset \{\beta_0,\dots,\beta_q\}$ and $\{\alpha_0,\dots,\alpha_p\}
\neq \{\beta_0,\dots,\beta_q\}$. Recall that by definition of the
category $\Delta^{ext}(\vec{D}^0)^{op}$, one necessarily has $\alpha_0
= \beta_0 = \widehat{0}$ and $\alpha_p = \beta_q = \widehat{1}$. Such
a simplex ${\underline{\beta} = (\beta_0,\dots,\beta_q)}$ can be 
written as an expression of the form 
\[\alpha_0.\underline{\delta_1}.\underline{\delta_2}\dots \underline{\delta_p}\] 
with $\alpha_{i}.\underline{\delta_{i+1}} \supseteqq
(\alpha_{i},\alpha_{i+1})$ for all $0\leq i\leq p-1$ and such that at
least for one $i$, one has $\alpha_{i}.\underline{\delta_{i+1}}
\supsetneqq (\alpha_{i},\alpha_{i+1})$. And since the small category
$\Delta^{ext}(\vec{D}^0)^{op}$ only contains commutative diagrams, one
obtains the homeomorphism \be \label{equfirst}
L_{(\alpha_0,\dots,\alpha_p)}\mathcal{F}\vec{D}\iso
\liminj_{\{\alpha_0,\dots,\alpha_p\} \subsetneqq
  \{\beta_0,\dots,\beta_q\}}\mathcal{F}\vec{D}_{(\beta_0,\dots,\beta_q)}.\ee

Let $\mathcal{E}$ be the set of subsets $S$ of $\{0,\dots,p-1\}$ such
that $S\neq\{0,\dots,p-1\}$. Let $I(S)$ be the full subcategory of
$\Delta^{ext}(\vec{D}^0)^{op}$ consisting of the objects
${\underline{\beta}} = (\beta_0,\dots,\beta_q)$ such that 
\begin{enumerate}
\item $\{\alpha_0,\dots,\alpha_p\}
\subsetneqq \{\beta_0,\dots,\beta_q\}$
\item for any $i\notin S$, one has $\alpha_{i}.\underline{\delta_{i+1}}
\supsetneqq (\alpha_{i},\alpha_{i+1})$.
\end{enumerate}

The full subcategory $\bigcup_{S\in \mathcal{E}} I(S)$ is exactly the
subcategory of $\Delta^{ext}(\vec{D}^0)^{op}$ consisting of the
objects ${\underline{\beta}}$ such that $\{\alpha_0,\dots,\alpha_p\}
\subsetneqq \{\beta_0,\dots,\beta_q\}$, that is to say the subcategory
calculating $L_{(\alpha_0,\dots,\alpha_p)} \mathcal{F}\vec{D}$. In
other terms, one has the isomorphism of spaces \be \label{equ1}
\liminj_{\bigcup_{S\in \mathcal{E}} I(S)}\mathcal{F}\vec{D}\iso
L_{\underline{\alpha}}\mathcal{F}\vec{D}.\ee

The full subcategory $I(S)$ of $\Delta^{ext}(\vec{D}^0)^{op}$ has a
final subcategory $\overline{I(S)}$ consisting of the
$\underline{\beta} = (\beta_0,\dots,\beta_q)$ such that 
\begin{enumerate}
\item $\{\alpha_0,\dots,\alpha_p\}
\subsetneqq \{\beta_0,\dots,\beta_q\}$
\item for any $i\notin S$, one has $\alpha_{i}.\underline{\delta_{i+1}}
\supsetneqq (\alpha_{i},\alpha_{i+1})$
\item for any $i \in S$, one has $\alpha_{i}.\underline{\delta_{i+1}}
= (\alpha_{i},\alpha_{i+1})$. 
\end{enumerate}
The subcategory $\overline{I(S)}$ is final in $I(S)$ because for any
object $\underline{\beta}$ of $I(S)$, there exists a unique
$\underline{\gamma}$ of $\overline{I(S)}$ and a unique arrow
$\underline{\beta} \longrightarrow \underline{\gamma}$.
Therefore, by Theorem~\ref{souslim}, one has the isomorphism
\be \label{equ2}
\liminj_{I(S)}\mathcal{F}\vec{D}\iso \liminj_{\overline{I(S)}}\mathcal{F}\vec{D}
\ee
since the comma category
$(\underline{\beta}\!\downarrow\!\overline{I(S)})$ is the one-object
category. For any object $\underline{\beta}$ of $\overline{I(S)}$, one has
\begin{align*}
& \mathcal{F}\vec{D}_{\underline{\beta}} & \\
&= \prod_{i=0}^{i=q-1} \P_{\beta_i,\beta_{i+1}}\vec{D} & \hbox{ by definition of $\mathcal{F}\vec{D}$}\\
&= \prod_{i=0}^{i=p-1} \mathcal{F}\vec{D}_{\alpha_{i}.\underline{\delta_{i+1}}}  & \hbox{ by definition of $\mathcal{F}\vec{D}$}\\
&= \lp \prod_{i \in S}\mathcal{F}\vec{D}_{(\alpha_{i},\alpha_{i+1})}\rp 
\p \lp \prod_{i \notin S} \mathcal{F}\vec{D}_{\alpha_{i}.\underline{\delta_{i+1}}}\rp & \hbox{ by definition of $S$.}\\
\end{align*}
Thus, since the category $\top$ is cartesian closed, one obtains 
{\small \begin{align*}
& \liminj_{\overline{I(S)}}\mathcal{F}\vec{D} & \\
& \iso \liminj_{\overline{I(S)}} \lp \lp \prod_{i \in S}\mathcal{F}\vec{D}_{(\alpha_{i},\alpha_{i+1})}\rp 
\p \lp \prod_{i \notin S} \mathcal{F}\vec{D}_{\alpha_{i}.\underline{\delta_{i+1}}}\rp\rp  & \\ 
&\iso \lp \prod_{i \in S}\mathcal{F}\vec{D}_{(\alpha_{i},\alpha_{i+1})}\rp \p \liminj_{\begin{array}{c}i \notin S \\ \alpha_{i}.\underline{\delta_{i+1}}
\supsetneqq (\alpha_{i},\alpha_{i+1})\end{array}}\lp \prod_{i \notin S} \mathcal{F}\vec{D}_{\alpha_{i}.\underline{\delta_{i+1}}}\rp & \\
&\iso \lp \prod_{i \in S}\mathcal{F}\vec{D}_{(\alpha_{i},\alpha_{i+1})}\rp \p \lp \prod_{i \notin S} \liminj_{\alpha_{i}.\underline{\delta_{i+1}}
\supsetneqq (\alpha_{i},\alpha_{i+1})} \mathcal{F}\vec{D}_{\alpha_{i}.\underline{\delta_{i+1}}}\rp & \hbox{ by Lemma~\ref{colimproduit}.}\\
\end{align*}}
Therefore, one obtains the isomorphism of topological spaces 
\be \label{equ3}
\liminj_{\overline{I(S)}}\mathcal{F}\vec{D}\iso 
\lp \prod_{i \in S}\mathcal{F}\vec{D}_{(\alpha_{i},\alpha_{i+1})}\rp 
\p  \lp \prod_{i \notin S} L_{(\alpha_{i},\alpha_{i+1})} \mathcal{F}\vec{D} \rp .
\ee
thanks to Isomorphism~(\ref{equfirst}).

If $S$ and $T$ are two elements of $\mathcal{E}$ such that $S\subset
T$, then there exists a canonical morphism of diagrams
$I(S)\longrightarrow I(T)$ inducing a canonical morphism of
topological spaces 
\[\liminj_{\underline{\beta}\in
I(S)}\mathcal{F}\vec{D}_{\underline{\beta}}
\longrightarrow  \liminj_{\underline{\beta}\in
I(T)}\mathcal{F}\vec{D}_{\underline{\beta}}.\] 
Therefore, by Equation~(\ref{equ2}) and Equation~(\ref{equ3}), the double colimit 
\[\liminj_{S\in \mathcal{E}}\lp \liminj_{I(S)} \mathcal{F}\vec{D}\rp\]
calculates the source of the morphism $i_{(\alpha_0,\alpha_1)}\square
\dots \square i_{(\alpha_{p-1},\alpha_p)}$ by
Theorem~\ref{calculpushout}.  It then suffices to prove the
isomorphism
\[\liminj_{S\in \mathcal{E}}\lp \liminj_{I(S)} \mathcal{F}\vec{D}\rp \iso  
\liminj_{\{\alpha_0,\dots,\alpha_p\}
\subsetneqq \{\beta_0,\dots,\beta_q\}}\mathcal{F}\vec{D}_{\underline{\beta}}\] 
to complete the proof. 
For that purpose, it suffices to construct two
canonical morphisms 
\[\liminj_{S\in \mathcal{E}}\lp \liminj_{I(S)} \mathcal{F}\vec{D}\rp \longrightarrow   
\liminj_{\{\alpha_0,\dots,\alpha_p\}
\subsetneqq \{\beta_0,\dots,\beta_q\}}\mathcal{F}\vec{D}_{\underline{\beta}}\] 
and 
\[\liminj_{\{\alpha_0,\dots,\alpha_p\}
\subsetneqq \{\beta_0,\dots,\beta_q\}}\mathcal{F}\vec{D}_{\underline{\beta}} 
\longrightarrow \liminj_{S\in \mathcal{E}}\lp \liminj_{I(S)} \mathcal{F}\vec{D}\rp .
\]
The first morphism comes from the isomorphism of
Equation~(\ref{equ1}). 
As for the second morphism, let us consider a
diagram of flows of the form:
\[
\xymatrix{
\mathcal{F}\vec{D}_{\underline{\beta}} \fd{} \fr{} & \liminj_{S\in \mathcal{E}}\lp \liminj_{I(S)} \mathcal{F}\vec{D}\rp\\
\mathcal{F}\vec{D}_{\underline{\gamma}} \ar@{->}[ru] &}
\]
One has to prove that it is commutative. Since one has $\bigcup_{S\in
\mathcal{E}} I(S)=\Delta^{ext}(\vec{D}^0)^{op}$, there exists $S\in \mathcal{E}$
such that $\underline{\gamma}$ is an object of $I(S)$. 
So
$\underline{\beta}$ is an object of $I(S)$ as well and there exists a
commutative diagram
\[
\xymatrix{
\mathcal{F}\vec{D}_{\underline{\beta}} \fd{} \fr{} &  \liminj_{I(S)} \mathcal{F}\vec{D}\\
\mathcal{F}\vec{D}_{\underline{\gamma}} \ar@{->}[ru] &}
\]
since the subcategory $\Delta^{ext}(\vec{D}^0)^{op}$ is
commutative. Hence the result. 
\epf

\section{Calculating the underlying homotopy type}
\label{calhoty}

\bth \label{cofibranthat}
Let $\vec{D}$ be a full directed ball.  Then the diagram of spaces
$\mathcal{F}{Q(\vec{D})}$ (where $Q$ is the cofibrant replacement
functor of $\dtop$) is Reedy cofibrant.
\eth

\bpf By Proposition~\ref{ll0} and since the
model category $\top$ is monoidal, one deduces that for any object
${\underline{\alpha}}$ of $\mathcal{F}{Q(\vec{D})}$, the topological
space $\mathcal{F}{Q(\vec{D})}_{\underline{\alpha}}$ is cofibrant.  By
Theorem~\ref{simplification} and by induction on the cardinal of the
set $\vec{D}^0$, it then suffices to prove that the continuous map
$L_{(\widehat{0},\widehat{1})}\mathcal{F}{Q(\vec{D})}
\longrightarrow \mathcal{F}{Q(\vec{D})}_{(\widehat{0},\widehat{1})}$
is a cofibration of topological spaces.

Let $X$ be an object of $\hda$ such that $X^0=\vec{D}^0$ and such that
the continuous map $L_{(\widehat{0},\widehat{1})}\mathcal{F}{X}
\longrightarrow \mathcal{F}{X}_{(\widehat{0},\widehat{1})}$ is a
cofibration of topological spaces. Consider a pushout diagram of flows
with $n\geq 0$ as follows:
\[
\xymatrix{
\glob(\mathbf{S}^{n-1})\fd{}\fr{\phi} & X\fd{}\\
\glob(\mathbf{D}^{n}) \fr{} & \cocartesien Y}
\] 
One wants to prove that the continuous map
$L_{(\widehat{0},\widehat{1})}\mathcal{F}{Y} \longrightarrow
\mathcal{F}{Y}_{(\widehat{0},\widehat{1})}$ is a cofibration of
topological spaces as well. One has the equality $X^0=Y^0$ since the
morphism $\glob(\mathbf{S}^{n-1})\longrightarrow\glob(\mathbf{D}^{n})$
restricts to the identity of $\{\widehat{0},\widehat{1}\}$ on the
$0$-skeletons and since the $0$-skeleton functor $X\mapsto X^0$
preserves colimits~\footnote{One has the canonical bijection $\set(X^0,Z)\iso\dtop(X,T(Z))$ where 
$T(Z)$ is the flow defined by $T(Z)^0=Z$ and for any
$(\alpha,\beta)\in Z\p Z$, $\P_{\alpha,\beta} T(Z)=\{0\}$.}. So one has
the commutative diagram
\[
\xymatrix{
L_{(\widehat{0},\widehat{1})}\mathcal{F}{X}\fr{}\ar@{^(->}[d] & L_{(\widehat{0},\widehat{1})}\mathcal{F}Y\fd{}\\
\mathcal{F}{X}_{(\widehat{0},\widehat{1})}\fr{}& \mathcal{F}Y_{(\widehat{0},\widehat{1})}.}
\]
There are two mutually exclusive cases:
\begin{enumerate}
\item $(\phi(\widehat{0}),\phi(\widehat{1}))=(\widehat{0},\widehat{1})$. 
One then has the situation
\[
\xymatrix{
L_{(\widehat{0},\widehat{1})}\mathcal{F}{X}\fr{=}\ar@{^(->}[d] & L_{(\widehat{0},\widehat{1})}\mathcal{F}Y\fd{}\\
\mathcal{F}{X}_{(\widehat{0},\widehat{1})}\ar@{^(->}[r]& \mathcal{F}Y_{(\widehat{0},\widehat{1})}}
\]
where the bottom horizontal arrow is a cofibration since it is a
pushout of the morphism of flows
$\glob(\mathbf{S}^{n-1})\longrightarrow\glob(\mathbf{D}^{n})$. So the
continuous map $L_{(\widehat{0},\widehat{1})}\mathcal{F}{Y} \longrightarrow
\mathcal{F}{Y}_{(\widehat{0},\widehat{1})}$ is a cofibration. 
\item
  $(\phi(\widehat{0}),\phi(\widehat{1}))\neq(\widehat{0},\widehat{1})$.
  Then, one has the pushout diagram of flows
\[
\xymatrix{
L_{(\widehat{0},\widehat{1})}\mathcal{F}{X}\fr{}\ar@{^(->}[d] & L_{(\widehat{0},\widehat{1})}\mathcal{F}Y\fd{}\\
\mathcal{F}{X}_{(\widehat{0},\widehat{1})}\fr{}& \cocartesien\mathcal{F}Y_{(\widehat{0},\widehat{1})}.}
\]
So the continuous map $L_{(\widehat{0},\widehat{1})}\mathcal{F}{Y}
\longrightarrow \mathcal{F}{Y}_{(\widehat{0},\widehat{1})}$ is again a
cofibration. In this situation, it may happen that
$L_{(\widehat{0},\widehat{1})}\mathcal{F}{X} =
L_{(\widehat{0},\widehat{1})}\mathcal{F}Y$.
\end{enumerate}
The proof is complete with Proposition~\ref{compgen}, and because the
canonical morphism of flows $\vec{D}^0 \longrightarrow \vec{D}$ is a
relative $I^{gl}$-cell complex, and at last because the property above
is clearly satisfied for $X=\vec{D}^0$.
\epf

\bth \label{cofibranthatcor}
Let $\vec{D}$ be a full directed ball.  Then the diagram of spaces
\[|\glob^{top}(\mathcal{F}{Q(\vec{D})})|\] (where $Q$ is the cofibrant
replacement functor of $\dtop$) is Reedy cofibrant.
\eth

\bpf 
The endofunctor of $\top$ defined by the mapping $Z \mapsto
|\glob^{top}(Z)|$ preserves colimits. Therefore, one has the
isomorphism
\[L_{(\widehat{0},\widehat{1})}|\glob^{top}(\mathcal{F}{Q(\vec{D})})| \iso 
\left|\glob^{top}\lp L_{(\widehat{0},\widehat{1})}\mathcal{F}{Q(\vec{D})} \rp \right|.\]
It remains to prove that this endofunctor preserves
cofibrations~\footnote{This functor is of course very close to the
pointed suspension functor. But it is not known how to view it as a
left adjoint, and therefore as a left Quillen functor.}. The proof
will be then complete thanks to Theorem~\ref{cofibranthat}.

The space $|\glob^{top}(Z)|$ is equal to the colimit of the diagram of
spaces $\mathcal{D}(Z)$ 
\[
\xymatrix{
\{0\}\p Z \ar@{->}[rd]\fd{} & & \{1\}\p Z \ar@{->}[ld] \fd{}\\
\{\widehat{0}\} & [0,1]\p Z & \{\widehat{1}\}.
}
\]
Let us consider the small category $\C$
\[
\xymatrix{
b \ar@{->}[rd]\fd{} & & d \ar@{->}[ld] \fd{}\\
a & c & e
}
\]
equipped with the Reedy structure 
\[
\xymatrix{
1 \ar@{->}[rd]\fd{} & & 1 \ar@{->}[ld] \fd{}\\
0 & 2 & 0.
}
\]
If $D$ is an object of the diagram category $\top^\C$, then the
latching spaces and the matching spaces of $D$ are equal to:
\begin{enumerate}
\item $L_a D= L_b D = L_d D = L_e D = \varnothing$
\item $L_c D= D_b \sqcup D_d$ 
\item $M_a D = M_e D = M_c D = \mathbf{1}$
\item $M_b D = D_a$ 
\item $M_d D = D_e$.
\end{enumerate}

A morphism of diagrams $D \longrightarrow E$ is a Reedy fibration if
\begin{enumerate}
\item $D_a \longrightarrow E_a \p_{M_a E} M_a D = E_a$ is a fibration
\item $D_e \longrightarrow E_e \p_{M_e E} M_e D = E_e$ is a fibration
\item $D_c \longrightarrow E_c \p_{M_c E} M_c D = E_c$ is a fibration
\item $D_b \longrightarrow E_b \p_{M_b E} M_b D = E_b \p_{E_a} D_a$ is a fibration
\item $D_d \longrightarrow E_d \p_{M_d E} M_d D = E_d \p_{E_e} D_e$ is a fibration. 
\end{enumerate}

Consider the categorical adjunction $\liminj :\top^\C \leftrightarrows
\top : \diag$. By the calculations above, if $X\longrightarrow Y$ is a
(resp. trivial) fibration of spaces, then $\diag(X) \longrightarrow
\diag(Y)$ is a (trivial) Reedy fibration. The colimit functor from
$\top^\C$ to $\top$ is therefore a left Quillen functor.

A morphism of diagrams $D \longrightarrow E$ is a Reedy cofibration if 
\begin{enumerate}
\item $D_a = D_a \sqcup_{L_a D} L_a E \longrightarrow E_a$ is a cofibration 
\item $D_b = D_b \sqcup_{L_b D} L_b E \longrightarrow E_b$ is a cofibration 
\item $D_d = D_d \sqcup_{L_d D} L_d E \longrightarrow E_d$ is a cofibration 
\item $D_e = D_e \sqcup_{L_e D} L_e E \longrightarrow E_e$ is a cofibration 
\item $D_c \sqcup_{(D_b \sqcup D_d)} (E_b\sqcup E_d) = D_c \sqcup_{L_c D} L_c E \longrightarrow E_c$ is a cofibration. 
\end{enumerate}

Now take a cofibration $Z_1 \longrightarrow Z_2$. Since the colimit
functor $\liminj:\top^\C \longrightarrow \top$ preserves cofibrations,
it then suffices to check that the morphism of diagrams
$\mathcal{D}(Z_1) \longrightarrow \mathcal{D}(Z_2)$ is a Reedy
cofibration. It then suffices to check the fifth condition above, that
is to say it suffices to prove that the continuous map
\[([0,1]\p Z_1) \sqcup_{\{0\}\p Z_1 \sqcup \{1\}\p Z_1} (\{0\}\p Z_2 \sqcup \{1\}\p Z_2) \longrightarrow ([0,1]\p Z_2) \] 
is a cofibration of topological spaces. It turns out that the latter
map is the pushout product (cf. Notation~\ref{notapushcal}) of the two
cofibrations $\{0,1\} \longrightarrow [0,1]$ and $Z_1 \longrightarrow
Z_2$. The proof is then complete because $\top$ is a monoidal model
category.
\epf

\bth \label{cofibrant2}
Let $\vec{D}$ be a full directed ball.  Then the diagram of flows
$\mathcal{G}{Q(\vec{D})}$ (where $Q$ is the cofibrant replacement
functor of $\dtop$) is Reedy cofibrant.
\eth

\bpf The argument is different from the one of
Theorem~\ref{cofibranthat}.  The flow $Q(\vec{D})$ is an object of
$\hda$. Therefore, the canonical morphism of flows
$\vec{D}^0\longrightarrow Q(\vec{D})$ is a priori a transfinite
composition of pushouts of elements of $I^{gl}_+ = I^{gl}\cup
\{R,C\}$. Since there is a bijection of sets $\vec{D}^0 \iso
Q(\vec{D})^0$, a pushout of $R:\{0,1\} \longrightarrow \{0\}$ or of
$C:\varnothing \longrightarrow \{0\}$ in the globular decomposition of
the relative $I^{gl}_+$-cell complex $\vec{D}^0\longrightarrow
Q(\vec{D})$ is necessarily without effect on $\vec{D}^0$. Thus, the
canonical morphism of flows $\vec{D}^0\longrightarrow Q(\vec{D})$ is a
transfinite composition of pushouts of elements of $I^{gl}$. So there
exists an ordinal $\lambda$ and a $\lambda$-sequence
$M:\lambda\longrightarrow\dtop$ such that $M_0=\vec{D}^0$,
$M_\lambda=Q(\vec{D})$ and for any $\mu<\lambda$, the morphism of
flows $M_\mu \longrightarrow M_{\mu+1}$ is a pushout of the inclusion
of flows $e_\mu:\glob(\mathbf{S}^{n_\mu-1})\subset
\glob(\mathbf{D}^{n_\mu})$ for some $n_\mu\geq 0$, that is one has the
pushout diagram of flows:
\[
\xymatrix{
\glob(\mathbf{S}^{n_\mu-1})\fd{}\fr{\phi_\mu}& M_\mu\fd{}\\
\glob(\mathbf{D}^{n_\mu})\fr{} & M_{\mu+1} .\cocartesien}
\]
Let $(\alpha_0,\dots,\alpha_p)$ be a simplex of
$\Delta^{ext}(\vec{D}^0)^{op}$. The relative $I^{gl}$-cell complex 
\[\vec{D}^0\longrightarrow
\mathcal{G}{Q(\vec{D})}_{(\alpha_0,\dots,\alpha_p)} =
\vec{D}\!\restriction_{[\alpha_0,\alpha_1]} * \dots *
\vec{D}\!\restriction_{[\alpha_{p-1},\alpha_p]} \] is a relative
$I^{gl}$-cell subcomplex which is the union of the globular cells
$e_\mu$ such that $[\phi_\mu(\widehat{0}),\phi_\mu(\widehat{1})]
\subset [\alpha_i,\alpha_{i+1}]$ for some $0\leq i<p$~\footnote{A
  $I^{gl}$-cell subcomplex is characterized by its cells since any
  morphism of $I^{gl}$ is an effective monomorphism of flows by
  \cite{model3} Theorem~10.6 and by \cite{ref_model2}
  Proposition~12.2.1.}.  So the subcomplex $\vec{D}^0\longrightarrow
\mathcal{G}{Q(\vec{D})}_{(\alpha_0,\dots,\alpha_p)}$ contains the
globular cells $e_\mu$ such that
$[\phi_\mu(\widehat{0}),\phi_\mu(\widehat{1})] \subset
[\alpha_0,\alpha_1]\sqcup \dots\sqcup [\alpha_{p-1},\alpha_p]$
($\sqcup$ meaning the disjoint union !).

We then deduce that all morphisms of the diagram
$\mathcal{G}{Q(\vec{D})}$ are inclusions of relative $I^{gl}$-cell
subcomplexes. Thus, the canonical morphism of flows
\[L_{(\alpha_0,\dots,\alpha_p)}\mathcal{G}{Q(\vec{D})} \longrightarrow
\mathcal{G}{Q(\vec{D})}_{(\alpha_0,\dots,\alpha_p)}\] is an inclusion
of relative $I^{gl}$-cell subcomplexes as well. More precisely, it is
equal to the transfinite composition of the inclusions of flows
$\glob(\mathbf{S}^{n_\mu-1})\subset \glob(\mathbf{D}^{n_\mu})$ such
that $[\phi_\mu(\widehat{0}),\phi_\mu(\widehat{1})] \subset
[\alpha_0,\alpha_1]\sqcup \dots\sqcup [\alpha_{p-1},\alpha_p]$ and such
that there does not exist any state $\alpha$ such that
$[\phi_\mu(\widehat{0}),\phi_\mu(\widehat{1})] \subset
[\alpha_0,\alpha_1] \sqcup \dots \sqcup [\alpha_{i},\alpha]\sqcup
[\alpha,\alpha_{i+1}] \sqcup \dots \sqcup [\alpha_{p-1},\alpha_p]$ and
$\alpha_i<\alpha<\alpha_{i+1}$.  \epf

The proof of Theorem~\ref{cofibrant2} also has the following
consequences:

\begin{cor} 
\label{topologisation}
Let $\vec{D}$ be a full directed ball.  Then there exists a diagram of
globular complexes
\[\mathcal{G}^{top}{Q(\vec{D})}:\Delta^{ext}(\vec{D}^0)^{op}\longrightarrow \gltop\] 
such that the composition by the functor $cat:\gltop \longrightarrow
\dtop$
\[\Delta^{ext}(\vec{D}^0)^{op}\longrightarrow \gltop \longrightarrow \dtop\]
is exactly the diagram $\mathcal{G}{Q(\vec{D})}$.  
\end{cor}

\bpf First of all, consider the flow $Q(\vec{D})$ and using
Theorem~\ref{remonte}, construct a globular complex $Q(\vec{D})^{top}$
such that $cat(Q(\vec{D})^{top})=Q(\vec{D})$. Let
$(\alpha_0,\dots,\alpha_p)$ be a simplex of
$\Delta^{ext}(\vec{D}^0)^{op}$. Then the globular complex
\[\mathcal{G}^{top}{Q(\vec{D})}_{(\alpha_0,\dots,\alpha_p)}\] is
defined as the globular subcomplex containing the globular cells of
$Q(\vec{D})^{top}$ such that the attaching map $\phi$ satisfies
$[\phi(\widehat{0}),\phi(\widehat{1})] \subset
[\alpha_0,\alpha_1]\sqcup \dots\sqcup [\alpha_{p-1},\alpha_p]$.  \epf

Let $(\alpha_0,\dots,\alpha_p)$ be a simplex of
$\Delta^{ext}(\vec{D}^0)^{op}$. The category of multipointed
topological spaces being cocomplete, one can consider the multipointed
topological space
\[L_{(\alpha_0,\dots,\alpha_p)}\mathcal{G}^{top}{Q(\vec{D})}.\] It
consists of the globular subcomplexes of $Q(\vec{D})^{top}$ containing
the globular cells such that the attaching map $\phi$ satisfies
$[\phi(\widehat{0}),\phi(\widehat{1})] \subset [\alpha_0,\alpha_1]\sqcup
\dots\sqcup [\alpha_{p-1},\alpha_p]$ and such that there exists a state
$\alpha$ such that $[\phi(\widehat{0}),\phi(\widehat{1})] \subset
[\alpha_0,\alpha_1] \sqcup \dots \sqcup [\alpha_{i},\alpha]\sqcup
[\alpha,\alpha_{i+1}] \sqcup \dots \sqcup [\alpha_{p-1},\alpha_p]$ and
$\alpha_i<\alpha<\alpha_{i+1}$. So the multipointed topological space
$L_{(\alpha_0,\dots,\alpha_p)}\mathcal{G}^{top}{Q(\vec{D})}$ is a
globular complex. And one obtains the equality
\[cat(L_{(\alpha_0,\dots,\alpha_p)} \mathcal{G}^{top}{Q(\vec{D})}) =
L_{(\alpha_0,\dots,\alpha_p)}\mathcal{G}{\vec{D}}.\]

\begin{cor} \label{cofibrant2cor} 
With the choices of Corollary~\ref{topologisation}. Let $\vec{D}$ be a
full directed ball.  Then the diagram of spaces
$|\mathcal{G}^{top}{Q(\vec{D})}|$ (where $Q$ is the cofibrant
replacement functor of $\dtop$) is Reedy cofibrant.
\end{cor}

\bpf Let $(\alpha_0,\dots,\alpha_p)$ be a simplex of
$\Delta^{ext}(\vec{D}^0)^{op}$. The continuous map 
\[|L_{(\alpha_0,\dots,\alpha_p)}
\mathcal{G}^{top}{Q(\vec{D})}| \longrightarrow |\mathcal{G}^{top}{Q(\vec{D})}_{(\alpha_0,\dots,\alpha_p)}|\] 
is a transfinite composition of pushouts of continuous maps of the
form \[|\glob^{top}(\mathbf{S}^{n-1})|\longrightarrow
|\glob^{top}(\mathbf{D}^{n})|\] with $n\geq 0$. The proof is complete
thanks to the proof of Theorem~\ref{cofibranthatcor}.
\epf

\bth \label{fin}
Let $\vec{D}$ be a full directed ball.  Then its underlying homotopy
type $|\vec{D}|$ is the one of the point.
\eth

\bpf We are going to make an induction on the cardinal of the poset
$\vec{D}^0$. If $\vec{D}^0=\{\widehat{0}<\widehat{1}\}$, then
$Q(\vec{D})=\glob(Z)$ for some topological space $Z$. By hypothesis,
the space $Z=\P_{\widehat{0},\widehat{1}}\vec{D}$ is contractible (and
cofibrant).  Therefore, the flows $\glob(Z)$ and $\glob(\{0\})$ are
S-homotopy equivalent. Thus, the globular complexes $\glob^{top}(Z)$
and $\glob^{top}(\{0\})$ are S-homotopy equivalent as well. Hence the
topological spaces $|\glob^{top}(Z)|$ and $|\glob^{top}(\{0\})|$ are
homotopy equivalent by Theorem~\ref{remonte}.  Now suppose that
$\vec{D}^0\backslash \{\widehat{0}<\widehat{1}\}$ is non-empty and
suppose the theorem proved for any full directed ball $\vec{E}$ such
that $\card(\vec{E}^0)<\card(\vec{D}^0)$.

By Proposition~\ref{pushmax} applied to the full directed ball
$Q(\vec{D})$, one has the pushout diagram of flows:
\[
\xymatrix{
\glob(L_{(\widehat{0},\widehat{1})}\mathcal{F}{Q(\vec{D})}) \fr{}\fd{} & L_{(\widehat{0},\widehat{1})}\mathcal{G}{Q(\vec{D})}\fd{}\\
\glob(\P_{\widehat{0},\widehat{1}}Q(\vec{D})) \fr{} & Q(\vec{D}). \cocartesien}
\]
One obtains the commutative diagram of globular complexes:
\[
\xymatrix{
\glob^{top}(L_{(\widehat{0},\widehat{1})}\mathcal{F}{Q(\vec{D})})\fr{}\fd{} & L_{(\widehat{0},\widehat{1})}\mathcal{G}^{top}{Q(\vec{D})}\fd{}\\
\glob^{top}(\P_{\widehat{0},\widehat{1}}Q(\vec{D})) \fr{} & \mathcal{G}^{top}Q(\vec{D}) \cocartesien}
\]
which must be a pushout of multipointed topological spaces by
Corollary~\ref{cofibrant2cor}. One can now pass to the underlying
topological spaces of all of these globular complexes and one obtains
the pushout diagram of topological spaces:
\[
\xymatrix{
L_{(\widehat{0},\widehat{1})}|\glob^{top}(\mathcal{F}{Q(\vec{D})})|=|\glob^{top}(L_{(\widehat{0},\widehat{1})}\mathcal{F}{Q(\vec{D})})|\fr{}\fd{} & L_{(\widehat{0},\widehat{1})}|\mathcal{G}^{top}{Q(\vec{D})}|\fd{}\\
|\glob^{top}(\P_{\widehat{0},\widehat{1}}Q(\vec{D}))| \fr{} & |\mathcal{G}^{top}Q(\vec{D})|(\widehat{0},\widehat{1}). \cocartesien}
\]
The top horizontal arrow is induced by the morphism of diagrams of
spaces \[|\glob^{top}(\mathcal{F}{Q(\vec{D})})| \longrightarrow
|\mathcal{G}^{top}{Q(\vec{D})}|.\] If we can prove that the top
horizontal arrow is a weak homotopy equivalence of topological spaces,
and since the continuous map
$L_{(\widehat{0},\widehat{1})}|\glob^{top}(\mathcal{F}{Q(\vec{D})})|\longrightarrow
|\glob^{top}(\P_{\widehat{0},\widehat{1}}Q(\vec{D}))| $ is a
cofibration of spaces by Theorem~\ref{cofibranthatcor}, then one will
be able to deduce the weak homotopy equivalence of spaces
$|\glob^{top}(\P_{\widehat{0},\widehat{1}}Q(\vec{D}))| \simeq
|\mathcal{G}^{top}Q(\vec{D})|(\widehat{0},\widehat{1})$ since the
model category $\top$ is left proper. Since the topological space
$|\glob^{top}(\P_{\widehat{0},\widehat{1}}Q(\vec{D}))|$ is
contractible, one will be then able to deduce that the space
$|\mathcal{G}^{top}Q(\vec{D})|(\widehat{0},\widehat{1})\simeq
|\vec{D}|$ is weakly contractible. And the proof will be finished.

The diagrams of topological spaces
$|\glob^{top}(\mathcal{F}{Q(\vec{D})})|$ and
$|\mathcal{G}^{top}{Q(\vec{D})}|$ are both Reedy cofibrant by
Theorem~\ref{cofibranthatcor} and Corollary~\ref{cofibrant2cor}. So
their restriction to the full subcategory
$\de(\Delta^{ext}(\vec{D}^0)^{op}_+\!\downarrow\! 
(\widehat{0},\widehat{1})) \iso \Delta^{ext}(\vec{D}^0)^{op}\backslash
\{(\widehat{0},\widehat{1})\}$ of $\Delta^{ext}(\vec{D}^0)^{op}$ is 
Reedy cofibrant as well. Thus, one obtains 
\begin{align*}
& L_{(\widehat{0},\widehat{1})}|\glob^{top}(\mathcal{F}{Q(\vec{D})})| & \\
&\iso 
\liminj_{\de(\Delta^{ext}(\vec{D}^0)^{op}_+\!\downarrow\! (\widehat{0},\widehat{1}))} |\glob^{top}(\mathcal{F}{Q(\vec{D})})| & \hbox{ by definition of the latching space}\\
&\simeq \holiminj_{\de(\Delta^{ext}(\vec{D}^0)^{op}_+\!\downarrow\! (\widehat{0},\widehat{1}))} |\glob^{top}(\mathcal{F}{Q(\vec{D})})|& \hbox{ by Corollary~\ref{calculbienholim} and by Theorem~\ref{cofibranthatcor}}\\\end{align*}
and 
\begin{align*}
& 
L_{(\widehat{0},\widehat{1})}|\mathcal{G}^{top}{Q(\vec{D})}| & \\
&\iso 
\liminj_{\de(\Delta^{ext}(\vec{D}^0)^{op}_+\!\downarrow\! (\widehat{0},\widehat{1}))}|\mathcal{G}^{top}{Q(\vec{D})}| & \hbox{ by definition of the latching space}\\
&\simeq \holiminj_{\de(\Delta^{ext}(\vec{D}^0)^{op}_+\!\downarrow\! (\widehat{0},\widehat{1}))}
|\mathcal{G}^{top}{Q(\vec{D})}| & \hbox{ by Corollary~\ref{calculbienholim} and by Corollary~\ref{cofibrant2cor}.} \\
\end{align*}

It then suffices to prove that for any simplex
$(\alpha_0,\dots,\alpha_p)$ of the latching category
$\de(\Delta^{ext}(\vec{D}^0)^{op}_+\!\downarrow\! 
(\widehat{0},\widehat{1}))$, the morphism of diagrams
\[|\glob^{top}(\mathcal{F}{Q(\vec{D})})|\longrightarrow
|\mathcal{G}^{top}{Q(\vec{D})}|\] induces a weak homotopy equivalence
\[|\glob^{top}(\mathcal{F}{Q(\vec{D})})|(\alpha_0,\dots,\alpha_p)\simeq 
|\mathcal{G}^{top}{Q(\vec{D})}|(\alpha_0,\dots,\alpha_p).\]

The topological space
$|\mathcal{G}^{top}{Q(\vec{D})}|_{(\alpha_0,\dots,\alpha_p)}$ is the
``concatenation'' 
\[|\mathcal{G}^{top}{Q(\vec{D})}|_{(\alpha_0,\alpha_1)}*\dots
*|\mathcal{G}^{top}{Q(\vec{D})}|_{(\alpha_{p-1},\alpha_p)}\] of $p$
topological spaces, that is where the final state of
$\mathcal{G}^{top}{Q(\vec{D})}_{(\alpha_i,\alpha_{i+1})}$ is
identified with the initial state of
$\mathcal{G}^{top}{Q(\vec{D})}_{(\alpha_{i+1},\alpha_{i+2})}$ for any
$i+2\leq p$.  The latter space is contractible by induction hypothesis
and since a finite join of well-pointed cofibrant contractible spaces
is contractible. The topological space
\[|\glob^{top}(\mathcal{F}{Q(\vec{D})})|_{(\alpha_0,\dots,\alpha_p)}\]
is contractible since the product of spaces
\[\P_{\alpha_0,\alpha_1}Q(\vec{D})
\p\dots \p \P_{\alpha_{p-1},\alpha_p}Q(\vec{D})\] 
is contractible since $\vec{D}$ is a full directed ball and since a
finite product of cofibrant contractible spaces is contractible.
\epf

The proof of Theorem~\ref{fin} implies the following theorem:

\begin{cor}
Let $\vec{D}$ be a loopless flow such that 
\begin{enumerate}
\item the poset $\vec{D}^0$ is finite and bounded with initial 
state $\widehat{0}$ and with final state  $\widehat{1}$
\item for any $(\alpha,\beta)\in \vec{D}^0$ such that $\alpha<\beta$ 
and $(\alpha,\beta) \neq (\widehat{0},\widehat{1})$, the topological
space $\P_{\alpha,\beta}\vec{D}$ is weakly contractible.
\end{enumerate}
Then the underlying homotopy type of $\vec{D}$ is homotopy equivalent
to the underlying homotopy type of
$\glob(\P_{\widehat{0},\widehat{1}}\vec{D})$: in other terms, one has 
$|\vec{D}| \simeq |\glob(\P_{\widehat{0},\widehat{1}}\vec{D})|$. 
\end{cor}

\section{Preservation of the underlying homotopy type}
\label{prehoty}

\bth \label{thth}
Let $f:X\longrightarrow Y$ be a generalized T-homotopy
equivalence. Then the morphism $|f|:|X|\longrightarrow |Y|$ 
is an isomorphism of $\ho(\top)$. 
\eth

\bpf First of all, let us suppose that $f$ is a pushout diagram of
flows of the form
\[
\xymatrix{
Q(F(P_1))\ar@{^{(}->}[d]_{Q(F(u))}\fr{}& X\fd{f}\\ Q(F(P_2)) \fr{} & Y
\cocartesien}
\]
where $P_1$ and $P_2$ are two finite bounded posets and where
$u:P_1\longrightarrow P_2$ belongs to $\mathcal{T}$. Let us factor the
morphism of flows $Q(F(P_1))\longrightarrow X$ as a composite of a
relative $I^{gl}_+$-cell complex $Q(F(P_1))\longrightarrow W$ followed
by a trivial fibration $W\longrightarrow X$.  Then one obtains the
commutative diagram of flows
\[
\xymatrix{
Q(F(P_1))\ar@{^{(}->}[d]_{Q(F(u))}\ar@{^{(}->}[r] & W \ar@{^{(}->}[d] \ar@{->>}[r]^{\simeq}& X \fd{f}\\
Q(F(P_2)) \ar@{^{(}->}[r] & \cocartesien T \fr{\simeq} & \cocartesien Y.}
\]
The morphism $T\longrightarrow Y$ of the diagram above is a weak
S-homotopy equivalence since the model category $\dtop$ is left proper
by \cite{2eme} Theorem~6.4. So the flows $W$ and $X$ (resp. $T$ and
$Y$) have the same underlying homotopy types by \cite{model2}
Proposition~VII.2.2 and we are reduced to the following situation:
\[
\xymatrix{
Q(F(P_1))\ar@{^{(}->}[d]_{Q(F(u))}\ar@{^{(}->}[r]& X\ar@{^{(}->}[d]_{f}\\
Q(F(P_2)) \ar@{^{(}->}[r] & Y. \cocartesien}
\]
The four morphisms of the diagram above are inclusions of
$I^{gl}_+$-cell complexes. So using the globular decompositions of the
flows $Q(F(P_1))$, $Q(F(P_2))$, $X$ and $Y$, there exist four
globular complexes $Q^{top}(F(P_1))$, $Q^{top}(F(P_2))$, $X^{top}$ and
$Y^{top}$ and a commutative diagram of globular complexes
\[
\xymatrix{
Q^{top}(F(P_1))\ar@{^{(}->}[d]_{}\ar@{^{(}->}[r]& X^{top}\ar@{^{(}->}[d]_{}\\
Q^{top}(F(P_2)) \ar@{^{(}->}[r] & Y^{top}. \cocartesien}
\]
which is a pushout diagram of multipointed spaces and whose image by
the functor $cat:\gltop\longrightarrow \dtop$ gives back the diagram
of flows above. Now by passing to the underlying topological spaces,
one obtains the pushout diagram of topological spaces
\[
\xymatrix{
|Q^{top}(F(P_1))|\ar@{^{(}->}[d]_{}\ar@{^{(}->}[r]& |X^{top}|\ar@{^{(}->}[d]_{}\\
|Q^{top}(F(P_2))| \ar@{^{(}->}[r] & |Y^{top}|. \cocartesien}
\]
The continuous map $|Q^{top}(F(P_1))|\longrightarrow
|Q^{top}(F(P_2))|$ is a trivial cofibration of topological spaces
since the morphism of posets $u: P_1 \longrightarrow P_2$ is
one-to-one. Thus, the continuous map $|X^{top}|\longrightarrow
|Y^{top}|$ is a trivial cofibration as well.

Now let us suppose that $f:X\longrightarrow Y$ is a transfinite
composition of morphisms as above. Then there exists an ordinal
$\lambda$ and a $\lambda$-sequence $Z:\lambda\longrightarrow \dtop$
with $Z_0=X$, $Z_\lambda=Y$ and the morphism $Z_0\longrightarrow
Z_\lambda$ is equal to $f$. Since for any $u\in\mathcal{T}$, the
morphism of flows $Q(F(u))$ is a cofibration, the morphism
$Z_\mu\longrightarrow Z_{\mu+1}$ is a cofibration for any
$\mu<\lambda$.  Since the model category $\dtop$ is left proper by
\cite{2eme} Theorem~6.4, there exists by \cite{ref_model2}
Proposition~17.9.4 a $\lambda$-sequence
$\widetilde{Z}:\lambda\longrightarrow \dtop$ and a morphism of
$\lambda$-sequences $\widetilde{Z}\longrightarrow Z$ such that for any
$\mu\leq \lambda$, the flow $\widetilde{Z}_\mu$ is an object of
$\hda$, such that each morphism $\widetilde{Z}_\mu
\longrightarrow\widetilde{Z}_{\mu+1}$ is a relative $I^{gl}_+$-cell
complex, and such that the morphism $\widetilde{Z}_\mu\longrightarrow
Z_\mu$ is a weak S-homotopy equivalence. Using the globular
decomposition of $\widetilde{Z}_0$, construct a globular complex
$\widetilde{Z}_0^{top}$ such that
$cat(\widetilde{Z}_0^{top})=\widetilde{Z}_0$. And by transfinite
induction on $\mu$, since each morphism $\widetilde{Z}_\mu
\longrightarrow\widetilde{Z}_{\mu+1}$ is a relative $I^{gl}_+$-cell
complex, construct a globular complex $\widetilde{Z}_\mu^{top}$ such
that $cat(\widetilde{Z}_\mu^{top})=\widetilde{Z}_\mu$. Then one
obtains a $\lambda$-sequence of topological spaces $\mu\mapsto
|\widetilde{Z}_\mu^{top}|$ whose colimit is the underlying topological
space of $\widetilde{Z}_\lambda^{top}$.

For any $\mu<\lambda$, the continuous map $|\widetilde{Z}_\mu^{top}|
\longrightarrow |\widetilde{Z}_{\mu+1}^{top}|$ is a trivial
cofibration of topological spaces. So the transfinite composition
$|\widetilde{Z}_0^{top}| \longrightarrow
|\widetilde{Z}_\lambda^{top}|$ is a trivial cofibration as well.

It remains the case where $f$ is a retract of a generalized
T-equivalence of the preceding kinds. The result follows from the fact
that everything is functorial and that the retract of a weak homotopy
equivalence is a weak homotopy equivalence. 
\epf

\section{Conclusion}

This new definition of T-homotopy equivalence seems to be well-behaved
because it preserves the underlying homotopy type of flows. For an
application of this new approach of T-homotopy, see the proof of an
analogue of Whitehead's theorem for the full dihomotopy relation in
\cite{hocont}.

\appendix

\section{Elementary remarks about flows}
\label{limelm}

This is a reminder of results of \cite{3eme}.

\bp \label{precision}
(\cite{model3} Proposition~15.1) If one has the pushout of flows
\[
\xymatrix{
\glob(\de Z) \fr{\phi}\fd{} & A \fd{} \\
\glob(Z) \fr{} & \cocartesien M }
\]
then the continuous map $\P A\longrightarrow \P M$ is a transfinite
composition of pushouts of continuous maps of the form
$\id\p\dots\p\id\p f \p\id \p\dots\p \id$ where
$f:\P_{\phi(\widehat{0}),\phi(\widehat{1})} A\longrightarrow T$ is the
canonical inclusion obtained with the pushout diagram of topological
spaces
\[
\xymatrix{
\de Z\fr{}\fd{} & \P_{\phi(\widehat{0}),\phi(\widehat{1})}A\fd{} \\
Z \fr{} & \cocartesien T. }
\]
\ep

\bp\label{ll05}
Let $Y$ be a flow such that $\P Y$ is a cofibrant topological
space. Let $f:Y\longrightarrow Z$ be a pushout of a morphism of
$I^{gl}_+$. Then the topological space $\P Z$ is cofibrant.
\ep

\bpf By hypothesis, $f$ is the pushout of a morphism of flows
$g\in I^{gl}_+$. So one has the pushout of flows
\[
\xymatrix{
A\fd{g} \fr{\phi} & Y\ar@{->}[d]^-{f}\\
B \ar@{->}[r]_-{\psi} & Z. \cocartesien
}
\]
If $f$ is a pushout of $C:\varnothing\subset \{0\}$, then $\P Z = \P
Y$. Therefore, the space $\P Z$ is cofibrant. If $f$ is a pushout of
$R:\{0,1\}\rightarrow \{0\}$ and if $\phi(0)=\phi(1)$, then $\P Z = \P
Y$ again.  Therefore, the space $\P Z$ is cofibrant again. If $f$ is a
pushout of $R:\{0,1\}\rightarrow \{0\}$ and if $\phi$ is one-to-one,
then one has the homeomorphism \beas \P Z&\iso&\P Y\sqcup
\bigsqcup_{r\geq 0} \lp
\P_{.,\phi(0)}Y\p\P_{\phi(1),\phi(0)}Y\p \P_{\phi(1),\phi(0)}Y\p \dots\hbox{($r$ times)}  \p \P_{\phi(1),.}Y  \rp \\
&&\sqcup \bigsqcup_{r\geq 0} \lp
\P_{.,\phi(1)}Y\p\P_{\phi(0),\phi(1)}Y\p \P_{\phi(0),\phi(1)}Y\p
\dots\hbox{($r$ times)} \p \P_{\phi(0),.}Y \rp.  \eeas Therefore, the
space $\P Z$ is again cofibrant since the model category $\top$ is monoidal.
It remains the case where $g$ is the inclusion
$\glob(\mathbf{S}^{n-1})\subset \glob(\mathbf{D}^n)$ for some $n\geq
0$. Consider the pushout of topological spaces
\[
\xymatrix{
\mathbf{S}^{n-1}\fd{g} \fr{\P \phi} & \P_{\phi(\widehat{0}),\phi(\widehat{1})}Y\ar@{->}[d]^-{f}\\
\mathbf{D}^n \ar@{->}[r]_-{\P \psi} & T. \cocartesien
}
\]
By Proposition~\ref{precision}, the continuous map $\P
Y\longrightarrow \P Z$ is a transfinite composition of pushouts of
continuous maps of the form $\id\p\id\p\dots \p f \p \dots \p
\id\p\id$ where $f$ is a cofibration and the identities maps are the
identity maps of cofibrant topological spaces. So it suffices to
notice that if $k$ is a cofibration and if $X$ is a cofibrant
topological space, then $\id_X\p k$ is still a cofibration since the
model category $\top$ is monoidal.  \epf

\bp\label{ll0} 
Let $X$ be a cofibrant flow.  Then for any $(\alpha,\beta)\in X^0\p
X^0$, the topological space $\P_{\alpha,\beta}X$ is cofibrant. 
\ep

\bpf 
A cofibrant flow $X$ is a retract of a $I^{gl}_+$-cell complex $Y$ and
$\P X$ becomes a retract of $\P Y$. So it suffices to show that $\P Y$
is cofibrant. Proposition~\ref{ll05} completes the proof.
\epf

\section{Calculating pushout products}
\label{calpush}

This is a reminder of results of \cite{3eme}.

\begin{lem} \label{colimproduit} 
Let $D:I\longrightarrow \top$ and $E:J\longrightarrow \top$ be two
diagrams in a complete cocomplete cartesian closed category. Let $D\p
E:I\p J:\longrightarrow \top$ be the diagram of topological spaces
defined by $(D\p E)(x,y):=D(x)\p E(y)$ if $(x,y)$ is either an object
or an arrow of the small category $I\p J$. Then one has $\liminj (D\p
E)\iso(\liminj D)\p (\liminj E)$.
\end{lem}

\bpf One has $\liminj (D\p E)\iso \liminj_i(\liminj_j D(i)\p E(j))$ by 
\cite{MR1712872}. And  one has $\liminj_j (D(i)\p E(j))\iso D(i)\p (\liminj E)$ 
since the category is cartesian closed. So $\liminj (D\p E)\iso
\liminj_i (D(i)\p (\liminj E))\iso (\liminj D)\p (\liminj E)$. \epf

\begin{nota} \label{notapushcal}
If $f:U\longrightarrow V$ and $g:W\longrightarrow X$ are two morphisms
of a complete cocomplete category, then let us denote by $f\square g:
(U\p X) \sqcup_{(U\p W)} (V\p W)\longrightarrow V\p X$ the {\rm
pushout product} of $f$ and $g$. The notation $f_0\square\dots \square
f_p$ is defined by induction on $p$ by $f_0\square\dots \square
f_p:=(f_0\square\dots \square f_{p-1})\square f_p$.
\end{nota}

\bth\label{calculpushout}  (Calculating a pushout product of several morphisms)
Let $f_i:A_i\longrightarrow B_i$ for $0\leq i\leq p$ be $p+1$ morphisms
of a complete cocomplete cartesian closed category $\C$. Let $S\subset
\{0,\dots,p\}$. Let 
\[C_p(S):=\lp \prod_{i\in S} B_i\rp \p \lp \prod_{i\notin S} A_i\rp.\]
If $S$ and $T$ are two subsets of $\{0,\dots,p\}$ such that $S\subset
T$, let $C_p(i_S^T):C_p(S)\longrightarrow C_p(T)$ be the morphism 
\[\lp\prod_{i\in S}\id_{B_i}\rp\p \lp\prod_{i\in T\backslash S} f_i\rp\p 
\lp\prod_{i\notin T} \id_{A_i}\rp. \]
Then:  
\begin{enumerate} 
\item the mappings $S\mapsto C_p(S)$ and $i_S^T\mapsto C_p(i_S^T)$ 
give rise to a functor from $\Delta(\{0,\dots,p\})$ (the order
complex of the poset $\{0,\dots,p\}$) to $\C$
\item there exists a canonical morphism 
{
\[\liminj_{{\begin{array}{c}S\subset \{0,\dots,p\}\\ S\neq  \{0,\dots,p\}\end{array}}} 
C_p(S)\longrightarrow C_p(\{0,\dots,p\}).\]}
and it is equal to the morphism $f_0\square\dots \square f_p$. 
\end{enumerate}
\eth

\bpf The first assertion is clear. Moreover, for any subset $S$ and $T$ of 
$\{0,\dots,p\}$ such that $S\subset T$, the diagram
\[
\xymatrix{
S\fr{}\fd{} & \{0,\dots,p\} \\
T\ar@{->}[ru]&}
\]
is commutative since there is at most one morphism between two objects
of the order complex $\Delta(\{0,\dots,p\})$, hence the existence of
the morphism {\[\liminj_{{\begin{array}{c}S\subset \{0,\dots,p\}\\
        S\neq \{0,\dots,p\}\end{array}}} C_p(S)\longrightarrow
  C(\{0,\dots,p\}).\]}

The second assertion is clear for $p=0$ and $p=1$. We are going to
prove it by induction on $p$. By definition, the morphism
$f_0\square\dots \square f_{p+1}$ is the canonical morphism from
{\tiny\[\lp\lp\liminj_{{\begin{array}{c}S\subset \{0,\dots,p\}\\ S\neq  \{0,\dots,p\}\end{array}}} 
C_p(S)\rp \p B_{p+1}\rp
\sqcup_{\lp\lp\liminj_{{\begin{array}{c}S\subset \{0,\dots,p\}\\
S\neq \{0,\dots,p\}\end{array}}} C_p(S)\rp \p A_{p+1}\rp} \lp
C_p(\{0,\dots,p\})\p A_{p+1}\rp\]} to $B_0\p\dots\p B_{p+1}$. Since the
underlying category is supposed to be cartesian closed, the functors
$M\mapsto M\p B_{p+1}$ and $M\mapsto M\p A_{p+1}$ both preserve
colimits. So the source of the morphism $f_0\square\dots \square f_{p+1}$ is equal to 
{\scriptsize\[\lp\liminj_{{\begin{array}{c}S\subset \{0,\dots,p\}\\ S\neq  \{0,\dots,p\}\end{array}}} 
(C_p(S) \p B_{p+1})\rp
\sqcup_{\lp\liminj_{{\begin{array}{c}S\subset \{0,\dots,p\}\\
S\neq \{0,\dots,p\}\end{array}}} (C_p(S) \p A_{p+1})\rp} \lp
C_p(\{0,\dots,p\})\p A_{p+1}\rp\]}
or in other terms to 
{\tiny\[\lp\liminj_{{\begin{array}{c}S\subset \{0,\dots,p\}\\ S\neq  \{0,\dots,p\}\end{array}}} 
C_{p+1}(S\cup\{p+1\})\rp
\sqcup_{\lp\liminj_{{\begin{array}{c}S\subset \{0,\dots,p\}\\
S\neq \{0,\dots,p\}\end{array}}} C_{p+1}(S)\rp} 
C_{p+1}(\{0,\dots,p\})\]}
or at last to 
{\[\lp\liminj_{{\begin{array}{c}S\subset \{0,\dots,p+1\}\\ S\neq  \{0,\dots,p\}\\p+1\in S \end{array}}} 
C_{p+1}(S)\rp
\sqcup_{\lp\liminj_{{\begin{array}{c}S\subset \{0,\dots,p\}\\
S\neq \{0,\dots,p\}\end{array}}} C_{p+1}(S)\rp} 
C_{p+1}(\{0,\dots,p\}).\]}

The notation $\de\Delta(\{0,\dots,p+1\})$ will represent the
simplicial order complex $\Delta(\{0,\dots,p+1\})$ with the simplex
$(0,\dots,p+1)$ removed.

Let us consider the small category $\mathcal{D}$
\[
\xymatrix{
1 && 3 \\
& 2\ar@{->}[lu]^{u}\ar@{->}[ru]_{v}&}
\]
and the composition of functors: 
\[\de\Delta(\{0,\dots,p+1\})\longrightarrow \mathcal{D}\longrightarrow *\]
where $*$ is the category with one object and one morphism and where
the functor $F:\de\Delta(\{0,\dots,p+1\})\longrightarrow \mathcal{D}$
is defined as follows:
\begin{enumerate}
\item The full subcategory of $\de\Delta(\{0,\dots,p+1\})$ of $S$ such that 
$S\subset \{0,\dots,p+1\}$, $S\neq \{0,\dots,p\}$ and $p+1\in S$ is
mapped to $1$ and the identity morphism $\id_1$ of $1$.
\item The full subcategory of $\de\Delta(\{0,\dots,p+1\})$ of $S$ such that 
$S\subset \{0,\dots,p\}$ and $S\neq \{0,\dots,p\}$ is mapped to $2$
and the identity morphism $\id_2$ of $2$.
\item $F(\{0,\dots,p\})=3$. 
\item Any morphism from $F^{-1}(2)$ to $F^{-1}(1)$ is mapped to $u$
\item Any morphism from $F^{-1}(2)$ to $F^{-1}(3)$ is mapped to $v$.
\end{enumerate}
The functor $F$ gives rise to the adjunction between diagram categories: 
\[\boxed{F_*:\C^{\de\Delta(\{0,\dots,p+1\})\backslash
\{\{0,\dots,p+1\}\}}\leftrightarrows \C^{\mathcal{D}}:F^*}\]
where $F^*(X)=X\circ F$. It is easily seen that its left adjoint $F_*$
(i.e. the left Kan extension) sends a diagram $X$ of
$\C^{\de\Delta(\{0,\dots,p+1\})}$ to the diagram:
\[
\xymatrix{
\liminj\limits_{{\begin{array}{c}S\subset \{0,\dots,p+1\}\\ S\neq  \{0,\dots,p\}\\p+1\in S \end{array}}} 
X(S) && X(\{0,\dots,p\}) \\
& \liminj\limits_{{\begin{array}{c}S\subset \{0,\dots,p\}\\
S\neq \{0,\dots,p\}\end{array}}} X(S)\ar@{->}[lu]^{}\ar@{->}[ru]_{}&}
\]
The functor $\mathcal{D}\longrightarrow *$ gives rise to the
adjunction
\[\boxed{\liminj_\mathcal{D}:\C^\mathcal{D}\leftrightarrows\ C:\diag_\mathcal{D}}\]
where $\diag_\mathcal{D}$ is the diagonal functor. By composition of the two adjunctions, 
one obtains the isomorphism
{\scriptsize\[\liminj_{{\begin{array}{c}S\subset \{0,\dots,p+1\}\\ S\neq  \{0,\dots,p+1\}\end{array}}} X\iso \lp\liminj_{{\begin{array}{c}S\subset \{0,\dots,p+1\}\\ S\neq  \{0,\dots,p\}\\p+1\in S \end{array}}} 
X(S)\rp
\sqcup_{\lp\liminj_{{\begin{array}{c}S\subset \{0,\dots,p\}\\
S\neq \{0,\dots,p\}\end{array}}} X(S)\rp} 
X(\{0,\dots,p\}).\]}
This completes the induction. \epf

\section{Mixed transfinite composition of pushouts and cofibrations}
\label{mixcomp}

This is a reminder of results of \cite{3eme}.

\bp \label{compgen}
Let $\mathcal{M}$ be a model category. Let $\lambda$ be an ordinal. 
Let $(f_\mu:A_\mu\longrightarrow B_\mu)_{\mu<\lambda}$ be a $\lambda$-sequence 
of morphisms of $\mathcal{M}$. Let us suppose that for any $\mu<\lambda$, the 
diagram of
objects of $\mathcal{M}$
\[
\xymatrix{
A_\mu \fr{} \fd{f_\mu} & A_{\mu+1}\fd{}\\
B_\mu \ar@{^(->}[r]  & B_{\mu+1}
}
\]
is either a pushout diagram, or $A_\mu\rightarrow A_{\mu+1}$ is an
isomorphism and such that for any $\mu<\lambda$, $B_\mu\longrightarrow
B_{\mu+1}$ is a cofibration. Then: if $f_0:A_0\longrightarrow B_0$ is
a cofibration, then $f_\lambda:A_\lambda\longrightarrow B_\lambda$ is
a cofibration as well, where of course $A_\lambda:=\liminj A_\mu$ and
$B_\lambda:=\liminj B_\mu$.  \ep

\bpf It is clear that if $f_\mu:A_\mu\longrightarrow B_\mu$ is a cofibration, then 
$f_{\mu+1}:A_{\mu+1}\longrightarrow B_{\mu+1}$ is a cofibration as
well. It then suffices to prove that if $\nu\leq\lambda$ is a limit
ordinal such that $f_\mu:A_\mu\longrightarrow B_\mu$ is a cofibration
for any $\mu<\nu$, then $f_\nu:A_\nu\longrightarrow B_\nu$ is a
cofibration as well.  Consider a commutative diagram
\[
\xymatrix{
A_\nu \fr{} \fd{f_\nu} & C\fd{}\\
B_\nu \fr{}  \ar@{-->}[ru]^{k}& D
}
\]
where $C\longrightarrow D$ is a trivial fibration of
$\mathcal{M}$. Then one has to find $k:B_\nu\longrightarrow C$ making
both triangles commutative. Recall that by hypothesis,
$f_\nu=\liminj_{\mu<\nu} f_\mu$. Since $f_0$ is a cofibration, there 
exists a map $k_0$ making both triangles of the diagram 
\[
\xymatrix{
A_0 \fr{} \fd{f_0} & C\fd{}\\
B_0 \fr{}  \ar@{-->}[ru]^{k_0}& D
}
\]
commutative. Let us suppose $k_\mu$ constructed. There are two
cases. Either the diagram
\[
\xymatrix{
A_\mu \fr{} \fd{f_\mu} & A_{\mu+1}\fd{}\\
B_\mu \ar@{^(->}[r]  & B_{\mu+1}
}
\] 
is a pushout, and one can construct a morphism $k_{\mu+1}$ making both
triangles of the diagram
\[
\xymatrix{
A_{\mu+1} \fr{} \fd{f_{\mu+1}} & C\fd{}\\
B_{\mu+1} \fr{}  \ar@{-->}[ru]^{k_{\mu+1}}& D
}
\]
commutative and such that the composite $B_\mu \longrightarrow
B_{\mu+1} \longrightarrow C$ is equal to $k_\mu$ by using the
universal property satisfied by the pushout. Or the morphism
$A_\mu\rightarrow A_{\mu+1}$ is an isomorphism. In that latter case,
consider the commutative diagram
\[
\xymatrix{
B_\mu \fr{} \ar@{^(->}[d] \fr{k_\mu} & C\fd{}\\
B_{\mu+1} \fr{}  & D
}
\] 
Since the morphism $B_\mu \longrightarrow B_{\mu+1}$ is a cofibration,
there exists $k_{\mu+1} : B_{\mu+1} \longrightarrow C$ making the two
triangles of the latter diagram commutative. So, once again, the
composite $B_\mu \longrightarrow B_{\mu+1} \longrightarrow C$ is equal
to $k_\mu$.

The map $k:=\liminj_{\mu<\nu} k_\mu$ is a solution.
\epf

%\bibliographystyle{alpha} 
%\bibliography{invariance2T}

\end{document}